\documentclass[12pt,a4paper]{article}
\newfont{\bbf}{msbm10 at 12pt}
\newfont\bbfsm{msbm10 at 9pt}

\def\N{{\mbox{\bbf N}}}

\def\Z{{\mbox{\bbf Z}}}
\def\R{{\mbox{\bbf R}}}
\def\C{{\mbox{\bbf C}}}

\def\Nsm{\mbox{\bbfsm N}}

\newfont{\script}{eusm10 at 12pt}

\def\ovl{\overline}
\def\phi1{\phi}
\def\phi{\varphi}
\def\eps{\varepsilon}
\def\Re{\mbox{\rm Re}}
\def\Im{\mbox{\rm Im}}
\def\id{\mbox{\rm id}}

\def\sm{\setminus}

\def\*{ {\mbox{\tt $\star$}} }
\def\0{ {\mbox{\tt 0}} }
\def\1{ {\mbox{\tt 1}} }
\def\2{ {\mbox{\tt 2}} }
\def\3{ {\mbox{\tt 3}} }
\def\4{ {\mbox{\tt 4}} }

\def\s{ {\mbox{\tt s}} }

\def\St#1#2 {{\small\tt\rule{0pt}{0pt}_{\mbox{#1}}^{\mbox{#2}} }}

\def\<{\prec}
\def\>{\succ}

\ifx\ThmNum\undefined
    \newtheorem{theorem}{Theorem}[section]
    
\else
    \newtheorem{theorem}{Theorem}
    
\fi

\newtheorem{proposition}[theorem]{Proposition}
\newtheorem{lemma}[theorem]{Lemma}
\newtheorem{definition}[theorem]{Definition}

\newtheorem{corollary}[theorem]{Corollary}

\def\proof{\par\medskip\noindent {\sc Proof. }}
\def\proofof #1 {\par\medskip\noindent {\sc Proof of #1. }}

\def\sketchof #1 {\par\medskip\noindent {\sc Sketch of proof of #1. }}
\def\Box{\framebox[10pt]{\rule{0pt}{3pt}}}
\def\nix{\rule{0pt}{2pt}}
\def\qed{\qedd\par\medskip\noindent}
\def\qedd{\nix\nolinebreak\hfill\hfill\nolinebreak$\Box$}
\def\remark{\par\medskip \noindent {\sc Remark. }}
\def\lineclear
     {\rule{0pt}{0pt}\nopagebreak\par\nopagebreak\noindent}

\newlength\captionwidth
\def\reminder #1 {{\sf #1}}
\def\hide #1 {}
\long\def\longhide #1 {}

\usepackage{amsmath,amsfonts,eucal,graphicx,color}
\usepackage{epsfig}

\def\H{{\mbox{\bbf H}}}

\def\theta{\vartheta}

\newfont\Cyr{wncyr10 at 12pt}

\def\Log{\mbox{Log}}

\def\n{\nonumber}
\def\S{\mathcal S}
\def\G{\mathcal G}
\def\D{\mathcal D}

\def\sm{\setminus}
\def\s{{\underline s}}

\def\Nsm{{\mbox{\bbfsm N}}}

\newtheorem{definitionlemma}[theorem]{Definition and Lemma}

\begin{document}

\title{{\Large
Parameter Rays for the Exponential Family}}
\author{Markus F\"orster, Dierk Schleicher}

\maketitle

\thispagestyle{empty}

\begin{abstract}
We investigate the set of parameters $\kappa\in\C$ for which the
singular orbit $(0,e^{\kappa},\ldots)$ of
$E_{\kappa}(z):=\exp(z+\kappa)$ converges to $\infty$. 
These parameters are organized in smooth curves in parameter space
called \emph{parameter rays}.
\end{abstract}

\tableofcontents

\section{Introduction}

This article is a contribution to extend the successful theory of iterated
polynomials to transcendental entire functions. Douady and
Hubbard~\cite{Orsay}, and many others since then, have shown that the dynamical
planes of polynomials can be studied in terms of {\em dynamic rays} and their
landing points. Similarly, the space of quadratic polynomials can be understood in
terms of the structure of the Mandelbrot set, which itself is studied in terms of
{\em parameter rays}. 

In recent years, it has become clear that dynamical rays also make sense for
transcendental entire functions. Specifically for exponential functions, they were
introduced in \cite{dk,dgh1}, and in \cite{sz} it was shown that dynamic rays
classify all {\em escaping points}: those points which converge to $\infty$ under
iteration of the map. These results are extended to larger classes of entire
functions in \cite{dt,ro,GuenterThesis}. 

\begin{figure}
\includegraphics[width=\textwidth]{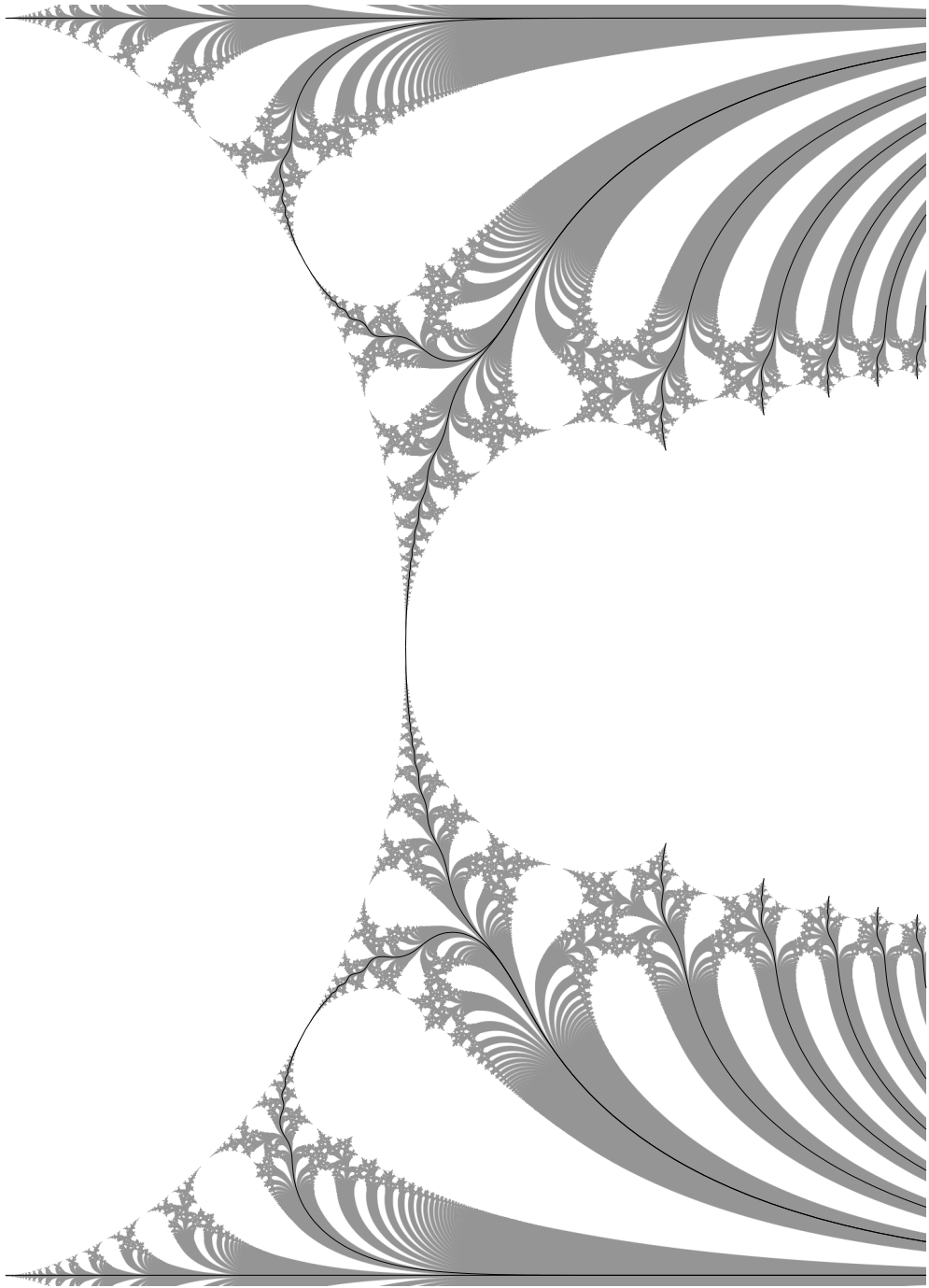}%
\caption{The parameter space of complex exponential maps $z\mapsto
e^z+\kappa$, with hyperbolic components in white, the bifurcation locus in
grey, and several parameter rays in black. (Picture courtesy of Lasse Rempe.) }
\label{FigParaRays}
\end{figure}

The simplest parameter space of transcendental functions is probably the space of
exponential functions; it has been studied in \cite{dgh2,el,br,rs1} and elsewhere
(see also \cite{fagella} and \cite{kk} for studies of other transcendental parameter
spaces). Similarly as for the Mandelbrot set, a systematic study of exponential
parameter space uses parameter rays and how parameter space is partitioned by
parameter rays landing at common points \cite{hab,rs1,rs2}; compare
Figure~\ref{FigParaRays}.

Parameter rays in the space of complex exponential maps are curves of parameters
for which the singular value ``escapes'', i.e.\ converges to $\infty$; such
parameters will be called ``escaping parameters''. Certain
parameter rays were constructed in \cite{dgh2} and more in \cite{hab}. In the
present paper, we construct and classify all parameter rays
(Theorem~\ref{fulllength}). In a sequel \cite{frs}, it will be shown how this helps
to classify all escaping parameters: every such parameter is either on a unique
parameter ray, or it is the landing point of a unique parameter ray with precisely
described combinatorics.

It turns out that the set of escaping parameters yields a nice dimension paradox:
the union of all parameter rays has Hausdorff dimension 1, while the set of only
those endpoints which are escaping parameters has dimension 2 \cite{bfs,frs}. This
is the parameter space analog to a well-known situation in the dynamical planes of
exponential maps \cite{k,sz}.

\begin{center}
{\bf Acknowledgements}
\end{center}
We would like to thank Lasse Rempe for many helpful discussions
and comments, as well as the dynamics workgroup at 
International University Bremen, especially Alex Kaffl, G\"unter
Rottenfu{\ss}er and Johannes R\"uckert. Furthermore, we would
like to thank the audiences in Paris, Warwick, and Oberwolfach for
interesting discussions, and the International University Bremen
for all the support.
Special thanks go to Niklas Beisert, who has contributed key ideas
to Section~\ref{sec_winnr}.

\section{Dynamic Rays}
\label{sec_dynrays}

\subsection{Notation and Definitions}
We investigate the family
\[
\{E_{\kappa}:\C\to\C\;,\ \;z\mapsto\exp(z+\kappa)\;\:|\
\kappa\in\C\}\;.
\]
Translating $\kappa$ by an integer multiple of $2\pi i$ yields the
same mapping, but slightly changes combinatorics. Therefore we
consider the complex plane as parameter space rather than the
cylinder $\C/2\pi i\Z$. The asymptotic value $0$ is the only
singular value of $E_\kappa$ (i.e. there are no other asymptotic
or critical values), and we call $\left(E_{\kappa}^{\circ
n}(0)\right)_{n\in\Nsm}=(0,e^\kappa,\exp(e^\kappa+\kappa),\dots)$
the \emph{singular orbit}.

As usual, the iterates of a function $f$ are denoted by
$f^{\circ(n+1)}(z):=f\circ f^{\circ n}(z)=f\circ\cdots\circ f(z)$,
with $f^{\circ 0}:=\id$. If $f$ is bijective, we define also
$f^{\circ (-n)}:=(f^{-1})^{\circ n}$. Let $\N:=\{0,1,2,3,\dots\}$,
$\N^*:=\N\setminus \{0\}$, $\C^*:=\C\setminus\{0\}$, and
$\C':=\C\setminus \R_0^{-}$. We shall write $B_r(z_0):=\{z\in\C:\
|z-z_0|<r\}$.

Let $\sigma:\S\to\S$ denote the shift map on the space
$\S:=\Z^{\Nsm^*}$ of sequences over the integers. We are going to
use
\[
F:\R^+_0\to \R^+_0\ ,\quad F(t):=e^t-1
\]
as a bijective model function for exponential growth.

The following discussion will take place in the dynamical plane of
a fixed map $E_\kappa$. 
We define
\hide{Recall the definitions}
\begin{eqnarray*}
I(E_{\kappa})&:=&\{z\in\C:\ |(E_{\kappa})^{\circ n}(z)|\to\infty
\mbox{ as }n\to\infty\}\;;\\ I&:=&\{\kappa\in\C:\ 0\in
I(E_\kappa)\}\;.
\end{eqnarray*}
\begin{lemma}[Characterization of Escaping Points]\label{esc_orb}
\lineclear For all $\kappa\in\C$,
\[
z\in I(E_\kappa)\quad\Longleftrightarrow\quad
\Re\left(E_{\kappa}^{\circ n}(z)\right) \to +\infty\quad \mbox{as
}n\to\infty\;.
\]
\end{lemma}
\begin{proof}
This follows from $|E_\kappa^{\circ(n+1)}(z)|=\exp(\Re
(E_{\kappa}^{\circ n}(z)+\kappa))$. \qed
\end{proof}

\begin{figure}
\begin{picture}(0,0)%
\includegraphics{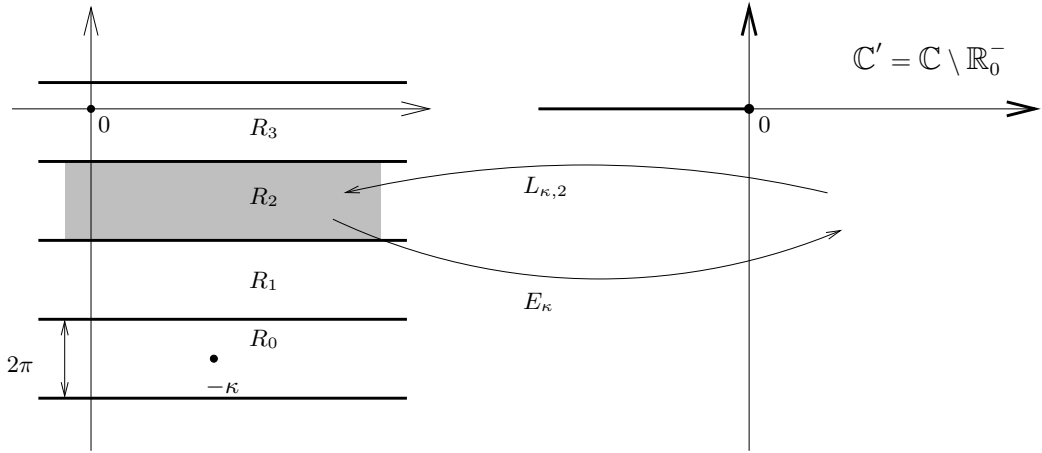}%
\end{picture}%
\setlength{\unitlength}{2901sp}%
\begingroup\makeatletter\ifx\SetFigFont\undefined%
\gdef\SetFigFont#1#2#3#4#5{%
  \reset@font\fontsize{#1}{#2pt}%
  \fontfamily{#3}\fontseries{#4}\fontshape{#5}%
  \selectfont}%
\fi\endgroup%
\begin{picture}(8842,3869)(1531,-4358)
\put(2311,-1606){\makebox(0,0)[lb]{\smash{\SetFigFont{9}{10.8}{\familydefault}{\mddefault}{\updefault}{\color[rgb]{0,0,0}$0$}%
}}}
\put(7951,-1606){\makebox(0,0)[lb]{\smash{\SetFigFont{9}{10.8}{\familydefault}{\mddefault}{\updefault}{\color[rgb]{0,0,0}$0$}%
}}}
\put(8761,-1066){\makebox(0,0)[lb]{\smash{\SetFigFont{10}{12.0}{\familydefault}{\mddefault}{\updefault}{\color[rgb]{0,0,0}$\C'=\C\sm\R^-_0$}%
}}}
\put(3601,-1636){\makebox(0,0)[lb]{\smash{\SetFigFont{9}{10.8}{\familydefault}{\mddefault}{\updefault}{\color[rgb]{0,0,0}$R_3$}%
}}}
\put(3601,-3436){\makebox(0,0)[lb]{\smash{\SetFigFont{9}{10.8}{\familydefault}{\mddefault}{\updefault}{\color[rgb]{0,0,0}$R_0$}%
}}}
\put(3601,-2221){\makebox(0,0)[lb]{\smash{\SetFigFont{9}{10.8}{\familydefault}{\mddefault}{\updefault}{\color[rgb]{0,0,0}$R_2$}%
}}}
\put(3601,-2941){\makebox(0,0)[lb]{\smash{\SetFigFont{9}{10.8}{\familydefault}{\mddefault}{\updefault}{\color[rgb]{0,0,0}$R_1$}%
}}}
\put(3231,-3841){\makebox(0,0)[lb]{\smash{\SetFigFont{9}{10.8}{\familydefault}{\mddefault}{\updefault}{\color[rgb]{0,0,0}$-\kappa$}%
}}}
\put(1531,-3661){\makebox(0,0)[lb]{\smash{\SetFigFont{9}{10.8}{\familydefault}{\mddefault}{\updefault}{\color[rgb]{0,0,0}$2\pi$}%
}}}
\put(5941,-2131){\makebox(0,0)[lb]{\smash{\SetFigFont{9}{10.8}{\familydefault}{\mddefault}{\updefault}{\color[rgb]{0,0,0}$L_{\kappa,2}$}%
}}}
\put(5941,-3121){\makebox(0,0)[lb]{\smash{\SetFigFont{9}{10.8}{\familydefault}{\mddefault}{\updefault}{\color[rgb]{0,0,0}$E_{\kappa}$}%
}}}
\end{picture}
\caption{The (static) partition and $L_{\kappa,j}$.}
\label{partit}
\end{figure}

To start with, we would like to endow the plane with dynamical
structure so as to obtain symbolic dynamics. On the slit plane
$\C'$, there is a biholomorphic branch $\Log:\C'\to \{z\in\C:|\Im
(z)|<\pi\}$ of the logarithm, which we will refer to as the
principal branch. Thus the branches of $E_\kappa^{-1}$ on $\C'$
are given by
\[
L_{\kappa,j}:\C' \to \C ;\quad L_{\kappa,j}(z)=\mbox{Log} z
-\kappa +2\pi ij\quad (j\in\Z)\;.
\]
Define the partition
\[
R_j := \{z\in \C: -\Im\kappa-\pi+2\pi j < \Im z <
-\Im\kappa+\pi+2\pi j\}\quad (j\in\Z)
\]
(see Figure \ref{partit}): these are the components of $\C\sm
E_{\kappa}^{-1}(\R^-)$. Clearly, every $L_{\kappa,j}$ maps $\C'$
biholomorphically onto $R_j$. The following definition gives rise
to symbolic dynamics and is the key idea for understanding the set
$I(E_\kappa)$.

\begin{definition}[External Addresses]\label{defextadr}
\lineclear Let $z=z_1\in \C$ be such that
$z_{n+1}:=E_\kappa^{\circ n}(z) \not \in \R^-$ for all $n\in\N$.
Then the \emph{external address} $\s(z)=(s_1,s_2,\dots)\in\S$ of
$z$ is defined to be the sequence of labels such that $z_n\in
 R_{s_n}$ for all $n \ge 1$.
\end{definition}

\remark The difference between a parameter $\kappa\in\C$ and its
translate $\kappa'=\kappa+2\pi i k$ (with $k\in\Z$) is a different
labeling: a point $z$ has external address $\s$ for the map
$E_{\kappa'}$ if and only if it has external address
$(s_1+k,s_2+k,s_3+k,\dots)$ for the map $E_{\kappa}$.

\begin{definition}[Exponential Boundedness]
\lineclear A sequence $\s \in\S$ is said to be \emph{exponentially
bounded} if there is an $x>0$, called \emph{growth parameter},
such that $|s_{k+1}|\le F^{\circ k}(x)$ for all $k\in\N$. The set
of exponentially bounded sequences is denoted by $\S_0$.
\end{definition}

The following lemma is taken from \cite{sz}, Lemma 2.4.

\begin{lemma}[External Addresses are Exponentially
Bounded]\label{extad_expbd} \lineclear For every $\kappa$ and
every $z\in\C$ there is an $x>0$ such that for all $n\ge 0$ we
have $|E_\kappa^{\circ n}(z)|\le F^{\circ n}(x)$. Thus every
sequence which is realized as an external address is exponentially
bounded. \qed
\end{lemma}

\subsection{Dynamic Rays}\label{secdynrays}

This section summarizes necessary results from \cite{sz} on
dynamic rays.

\begin{definitionlemma}[Minimal Potential]\label{minpot}
\lineclear For every $\s=(s_1,s_2,\dots) \in \S$ define the
\emph{minimal potential} of $\s$ by
\[
t_{\s}:=\inf\left\{ x>0:\; \limsup_{n \to \infty} \frac
{|s_n|}{F^{\circ (n-1)}(x)} = 0 \right\}\;,
\]
and furthermore define
\[
t_{\s}^*:=\sup_{n\ge 1} F^{\circ (-n+1)}(|s_n|)\;.
\]
These definitions lead to the following properties:
\begin{itemize}
\item $t_{\s}=\limsup_{n\ge 1} F^{\circ (-n+1)}(|s_n|)\le
t_{\s}^*$; \item If $\s\in\S_0$, then $t_{\s}\le t_{\s}^*<
\infty$; otherwise $t_{\s}=t_{\s}^*=\infty$; \item $|s_{n+1}|\le
F^{\circ n}(t_{\s}^*)$ for all $n\in\N$; \item $\forall t>t_{\s}\
\exists N_0\in\N\ \forall N\ge N_0:\quad F^{\circ N}(t)>
t_{\sigma^N\s}^*\;.$ In particular, for every $t>t_{\s}$ there
exists $N_0\in\N$ such that $\forall N\ge N_0:\ |s_{N+1}|\le
F^{\circ N}(t)$.
\end{itemize}
\end{definitionlemma}

\begin{proof}
Let $t_{\s}':=\limsup_{n\ge 1} F^{\circ (-n+1)}(|s_n|)$ and
$L(x):=\limsup_{n\ge 1}\frac{|s_n|}{F^{\circ n}(x)}$. Every
$x<t_{\s}$ satisfies $L(x)=\infty$, but $L(t_{\s}')= 1$. That
shows $t_{\s}'\ge t_{\s}$. Similarly, every $x$ such that $L(x)=0$
satisfies $x\ge t_{\s}'$, and thus $t_{\s}'\le t_{\s}$.

The second and third items follow directly from the definitions.
For the last item observe that $t_{\sigma^N\s}^*=F^{\circ
N}\left(\sup_{n\ge N+1} F^{\circ (-n+1)}(|s_{n}|)\right)$. Since
$t>\limsup_{n\ge 1} F^{\circ (-n+1)}(|s_n|)$, there is an $N_0$
such that for all $N\ge N_0$ we have $t>\sup_{n\ge N+1}F^{\circ
(-n+1)}(|s_n|)$. \qed
\end{proof}

According to Definition 6.7 in \cite{sz}, we divide $\S_0$ into
so-called \emph{slow} and \emph{fast} sequences: a sequence
$\s\in\S_0$ is called slow if it has a growth parameter $x$ which
works for infinitely many shifts $\sigma^n\s$ of $\s$ as well;
otherwise $\s$ is called fast. Consider the set $X\subset
\S_0\times\R_0^+$ defined by
\[
X:=\{(\s,t)\in\S_0\times\R:\ t>t_{\s}\}\cup
\bigcup_{\s\in\S_0\;\mbox{\scriptsize fast}} \{(\s,t_{\s})\}\;.
\]
Let us endow $X$ with the product topology induced by the discrete
topology on $\S_0$ and the standard topology on $\R$. (See
\cite{r2} for a deeper discussion of the topology of the sets
$I(E_{\kappa})$.)

\begin{theorem}[Dynamic Rays]\label{dynrays}
\lineclear \begin{enumerate} \item For every $\kappa\not\in I$
there is a continuous bijection
\[
g^{\kappa}:\;X\to I(E_{\kappa})\;.
\]
If $\kappa\in I$ then we have to restrict the map $g^{\kappa}$:
there is a preferred pair $(\s^{\kappa},t^{\kappa})$ and a set
\[
X^{\kappa}:=X\setminus\{(\s,t):\;\exists n\ge 1:
\sigma^n\s=\s^{\kappa}\mbox{ and }F^{\circ n}(t)\le t^{\kappa}\}
\]
with a continuous injection
$g^{\kappa}:X^{\kappa}\to I(E_{\kappa})$ with
$g^{\kappa}(\s^{\kappa},t^{\kappa})=0$. For every $z\in
I(E_{\kappa})\setminus g^{\kappa}(X^{\kappa})$ there is a unique
$n\ge 1$ and a $t<t^{\kappa}$ such that $E_{\kappa}^{\circ
n}(z)=g^{\kappa}(\s^{\kappa},t)$.

\item For every $\s\in\S_0$ and $\kappa\in\C$, we define the curve
\[g_{\s}^{\kappa}(t):=g^{\kappa}(\s,t)\qquad (t>t_{\s})\]
wherever this is defined. For fixed $t>t_{\s}$,
$g_{\s}^{\kappa}(t)$ depends analytically on $\kappa$ (on
$\kappa$-disks where $g_{\s}^{\kappa}(t)$ is defined). For fixed
$\kappa$ the curves have the following properties.
\begin{enumerate}
\item For all $t,\s$ where $g_{\s}^{\kappa}(t)$ is defined we have
\[
E_\kappa \circ g_{\s}^{\kappa}(t) = g_{\sigma\s}^{\kappa}\circ
F(t)\;.
\]
\item \label{defined}Suppose $K\ge |\kappa|$. Then
$g_{\s}^{\kappa}(t)$ is defined for all $t\ge
t_{\s}^K:=t_{\s}^*+2\log(K+3)$, and
\begin{eqnarray}
g_{\s}^{\kappa}(t) &=& t-\kappa+2\pi is_1 +r_{\kappa,\s}(t) \n\\
\mbox{with}\quad |r_{\kappa,\s}(t)| &\le&
2e^{-t}(K+2\pi|s_2|+12)<5\;.\label{rK_est}
\end{eqnarray}
\item The orbit $z_n:=E_{\kappa}^{\circ (n-1)}(z_1)$ of any
$z_1=g_{\s}^{\kappa}(t)$ satisfies
\begin{equation}\label{lowpotasymp}
z_n=F^{\circ (n-1)}(t)-\kappa+2\pi is_{n} +O\left(F^{\circ
(-n)}(t)\right)\quad \mbox{as } n\to\infty\;.
\end{equation}
\end{enumerate}
\end{enumerate}
\end{theorem}

\begin{proof}
Everything can be found in \cite{sz} (Prop.~3.2, Theorem 4.2, and
Prop.~4.4), except for the second inequality in (\ref{rK_est}): we
get for all $t\ge t_{\s}^K$ (using $x=t_{\s}^*$, $A=1$, $C<1.5$
for the variables appearing in \cite{sz})
\begin{eqnarray*}
|r_{\kappa,\s}(t)|\le \frac{2K+4+4\pi
F(t_{\s}^*)+6\pi}{e^{t_{\s}^*}(K+3)^2}\le\frac{2K+4\pi+24}{6K+9}<5\;.
\end{eqnarray*}

\vspace{-1cm} \qed
\end{proof}

\remark We call the curves $g_{\s}^{\kappa}(t)$ ($t>t_{\s}$)
\emph{dynamic rays} at external address $\s$. Viana \cite{v}
showed that they are $C^{\infty}$-smooth. (For $C^2$, see Section
\ref{sec_bound}.) For a fast sequence $\s$, the point
$g^{\kappa}(\s,t_{\s})$ (if defined) is called the \emph{endpoint}
of the dynamic ray $g_{\s}^{\kappa}$. See \cite{r2} for a
discussion of smoothness in the endpoints. The variable $t$ is
referred to as the \emph{potential} (in analogy to the terminology
for polynomial external rays).

\remark By Definition and Lemma \ref{minpot}, for every
$z=g_{\s}^{\kappa}(t)$ there is an $N$ such that $F^{\circ
(N-1)}(t)\ge t_{\sigma^N\s}^{|\kappa|}$, i.e.~for all $n\ge N$,
$E_{\kappa}^{\circ n}(z)=g_{\sigma^n\s}^{\kappa}(F^{\circ
(n-1)}(t))$ satisfies the condition required for (\ref{rK_est}). 
This is crucial for a lot of arguments given in this paper
(see Definition~\ref{raytails} and Lemma~\ref{Lem:FromTailsToRays}).

\section{Construction of Parameter Rays}

Now we want to turn our attention to \emph{parameter space}, the
set of parameters $\kappa$. We are interested in the set $I$ of
\emph{escaping parameters}, which are those parameters for which
the singular value $0$ is an escaping point. Again, this
investigation will lead to curves, called \emph{parameter rays},
which parameterize the escaping parameters by the external address
$\s$ and the potential $t$ (i.e. speed of escape) of the singular
orbit under $E_\kappa$.

Based on dynamic rays, we start the construction of the parameter
rays at large potentials, where it is comparably easy to find an
escaping parameter $\kappa$ with given combinatorial data
$(\s,t)$. Then we will extend these \emph{parameter ray tails}
onto the full domain $(t_{\s},\infty)$ of potentials.

\begin{definition}[Pointwise Definition of Parameter Rays]\label{defpara}
\lineclear For an external address $\s\in\S_0$ and a potential
$t>t_{\s}$, let
\[
\D_{\s}(t):=\C\setminus \{ \kappa\in I:\; \exists n\ge 1,\;
t^{\kappa}\ge F^{\circ n}(t):\
 g_{\sigma^n\s}^\kappa(t^{\kappa})=0\}
\]
be the set of parameters for which $g_{\s}^{\kappa}(t)$ is defined
in the sense of Theorem \ref{dynrays}, and
\begin{equation}\label{eqpara}
\G_{\s}(t):=\{\kappa\in\D_{\s}(t):\; g_{\s}^{\kappa}(t)=0\}\ .
\end{equation}
\end{definition}

\remark As our main result (Theorem \ref{fulllength}), we will
show that for every exponentially bounded external address
$\s\in\S_0$ and for every potential $t>t_{\s}$ we have
$|\G_{\s}(t)|=1$, and that the unique map
$G_{\s}:(t_{\s},\infty)\to I$, $t\mapsto \kappa\in \G_{\s}(t)$ is
a curve. This curve will be called the \emph{parameter ray} at
external address $\s$.

\subsection{Parameter Ray Tails}

\begin {proposition}[Existence of Parameter Ray Tails]\label{parrayend}
\lineclear For every $\s \in \S_0$ there is a constant
$T_{\s}>t_{\s}$ and a unique map $G_{\s}:[T_{\s},\infty) \to \C$,
called \emph{parameter ray tail}, such that for all $t\ge T_{\s}$
\[G_{\s}(t)\in \G_{\s}(t)\quad\mbox{and}\quad |G_{\s}(t)|<2\pi t\;.
\]
The parameter $G_{\s}(t)$ is a simple root of the map
$\kappa\mapsto g_{\s}^{\kappa}(t)$, and the parameter ray tails
carry the asymptotics
\[
G_{\s}(t)=t+ 2\pi i s_1 +R_{\s}(t)
\]
with $|R_{\s}(t)| < 2e^{-t}(2\pi t+2\pi |s_2|+12)<5$.
\end{proposition}

\begin{proof}
Define $T_{\s}:=20+2t_{\s}^*$. Consider an arbitrary fixed
potential $t\ge T_{\s}$ and define $K:=2\pi t$. We will show that
the map $h:\kappa\mapsto g_{\s}^{\kappa}(t):B_K(0)\to \C$ is
well-defined on the open disk $B_K(0)$ and contains a unique root
$\kappa_0$.

\begin{figure}
\begin{picture}(0,0)%
\includegraphics{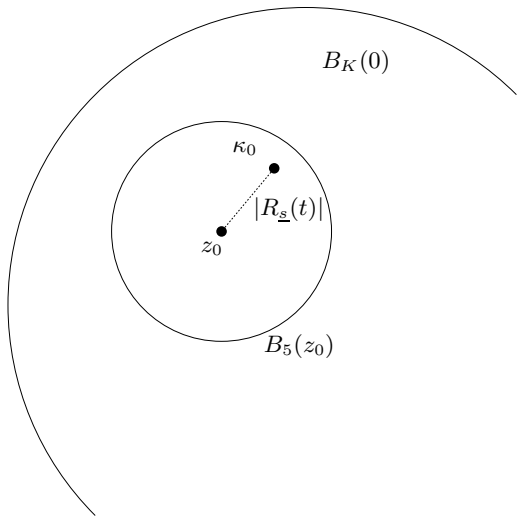}%
\end{picture}%
\setlength{\unitlength}{2901sp}%
\begingroup\makeatletter\ifx\SetFigFont\undefined%
\gdef\SetFigFont#1#2#3#4#5{%
  \reset@font\fontsize{#1}{#2pt}%
  \fontfamily{#3}\fontseries{#4}\fontshape{#5}%
  \selectfont}%
\fi\endgroup%
\begin{picture}(4363,4363)(146,-4119)
\put(1801,-1861){\makebox(0,0)[lb]{\smash{\SetFigFont{9}{10.8}{\familydefault}{\mddefault}{\updefault}{\color[rgb]{0,0,0}$z_0$}%
}}}
\put(2836,-286){\makebox(0,0)[lb]{\smash{\SetFigFont{9}{10.8}{\familydefault}{\mddefault}{\updefault}{\color[rgb]{0,0,0}$B_K(0)$}%
}}}
\put(2341,-2716){\makebox(0,0)[lb]{\smash{\SetFigFont{9}{10.8}{\familydefault}{\mddefault}{\updefault}{\color[rgb]{0,0,0}$B_{5}(z_0)$}%
}}}
\put(2071,-1006){\makebox(0,0)[lb]{\smash{\SetFigFont{9}{10.8}{\familydefault}{\mddefault}{\updefault}{\color[rgb]{0,0,0}$\kappa_0$}%
}}}
\put(2251,-1546){\makebox(0,0)[lb]{\smash{\SetFigFont{9}{10.8}{\familydefault}{\mddefault}{\updefault}{\color[rgb]{0,0,0}$|R_{\s}(t)|$}%
}}}
\end{picture}
\caption{The setting in the proof of Proposition \ref{parrayend}.}
\label{Gproof}
\end{figure}

Indeed, by Theorem \ref{dynrays} (\ref{defined}), this map is
well-defined, because $t\ge T_{\s}\ge 20$ implies $t/2>2\log(2\pi
t+3)$ and thus
\begin{eqnarray*}
t &=& t/2+t/2\;>\;2\log(2\pi t+3)+t_{\s}^*= t_{\s}^K\;.
\end{eqnarray*}
In fact, this implies by (\ref{rK_est}) that for all $\kappa\in
B_K(0)$
\begin{equation}\label{gskappa_est}
g_{\s}^{\kappa}(t)=t-\kappa+2\pi is_1 + r_{\kappa,\s}(t)
\quad\mbox{with } |r_{\kappa,\s}(t)|<5\;.
\end{equation}
Now for given $\kappa$, define $z_0:=t+ 2\pi i s_1$, so that
$g_{\s}^{\kappa}(t)=z_0-\kappa+r_{\kappa,\s}(t)$. Since
$|r_{\kappa,\s}(t)|<5$, we have $g_{\s}^{\kappa}(t) \neq 0$ for
$|z_0-\kappa|\geq 5$. Within the disk $B_K(0)$, the only
parameters $\kappa$ with $g_{\s}^{\kappa}(t)=0$ are thus contained
in the disk $B_{5}(z_0)$.

Note that $B_6(z_0)\subset B_K(0)$, because every $\kappa\in
B_6(z_0)$ satisfies $$|\kappa|\leq t+2\pi|s_1|+6 \le t+2\pi
t_{\s}^* + 6 \le t+\pi T_{\s}+6<2\pi t=K\;.$$ By
(\ref{gskappa_est}), $h(\partial B_{5.5}(z_0))$ winds exactly once
around $0$. Analyticity of $h$ and Rouch\'{e}'s theorem imply
therefore that there is exactly one $\kappa_0=:G_{\s}(t)$
(counting multiplicities) with $|\kappa_0|<K$ for which
$h(\kappa_0)=0$. This is a simple root of $\kappa\mapsto
g_{\s}^{\kappa}(t)$.

Since $g_{\s}^{\kappa_0}(t)=t-\kappa_0+2\pi
is_1+r_{\kappa_0,\s}(t)=0$, Theorem \ref{dynrays} (\ref{defined})
yields $|G_{\s}(t)-t-2\pi i s_1|\le |r_{2\pi t,\s}(t)|$. \qed
\end{proof}

Notice that we have not yet shown that $G_{\s}$ is a curve. The
following proposition implies uniqueness of the parameter ray
tails (without the restriction on $|\kappa|$) and is the main
argument for extending these ray tails onto the full domain of
definition $(t_{\s},\infty)$. We defer the proof, which is the
technical heart of this paper, to Section \ref{sec_bound}.

\begin{proposition}[A Bound on the Growth of Parameter Rays]\label{parabound}
\lineclear For every $\s\in\S_0$ there is a continuous function
$\xi_{\s}:(t_{\s},\infty)\to\R$ such that for every $t>t_{\s}$
\[
\G_{\s}(t)\subset B_{\xi_{\s}(t)}(0)\ .
\]
Moreover, for sufficiently large $t$ we can choose
$\xi_{\s}(t)=2t$.
\end{proposition}

\begin{corollary}[Parameter Ray Tails Are Unique]\label{parendunique}
\lineclear For every $\s \in\S_0$ and every sufficiently large
$t$, $\G_{\s}(t)$ contains exactly one element $\kappa_0$. The
multiplicity of $\kappa_0$ as a root of $\kappa\mapsto
g_{\s}^{\kappa}(t)$ is $1$.
\end{corollary}

\begin{proof}
By Proposition \ref{parabound} we have $\G_{\s}(t)\subset
B_{2t}(0)\subset B_{2\pi t}(0)$, and the claim follows from
Proposition \ref{parrayend}. \qed
\end{proof}

\subsection{Parameter Rays at Their Full Length}

\begin{lemma}[The Domain of Definition of $\kappa\mapsto
g_{\s}^{\kappa}(t)$]\label{defdomain} \lineclear Let $\s\in\S_0$
be an external address.
\begin{enumerate}
\item \label{existsN} For every bounded set $\Lambda$ of
parameters and every compact interval $J\subset (t_{\s},\infty)$
of potentials there is an $N\in\N$ such that
$$ \forall n\ge N,\;\forall\kappa\in \Lambda,\;\forall t\in J:\quad
\kappa\in \D_{\sigma^n\s}(F^{\circ n}(t))\;.$$

\item Let $t_0>t_{\s}$. For every $\kappa_0\in \D_{\s}(t_0)$ there
are neighborhoods $J\subset\R$ and $\Lambda\subset\C$ of $t_0$ and
$\kappa_0$ respectively such that
$$ \forall t\in J,\;\forall \kappa\in\Lambda:\quad \kappa\in\D_{\s}(t)\;.$$
In particular, $\D_{\s}(t)$ is open for every $t>t_{\s}$.
\end{enumerate}
\end{lemma}

\begin{proof}
Recall from Theorem \ref{dynrays} that if $t>t_{\s}$, the only
possible reason for $g_{\s}^{\kappa}(t)$ not to be defined is the
existence of an $n\ge 1$ such that
$g_{\sigma^n\s}^{\kappa}(t_0)=0$ with $t_0\ge F^{\circ n}(t)$.

For the first claim, let $K:=\sup _{\kappa\in \Lambda}|\kappa|$.
Take $N$ big enough such that $F^{\circ
N}(t)>\max\{t_{\sigma^N\s}^K,K+5\}$ for all $t\in J$. By Theorem
\ref{dynrays}, this implies for all $\kappa\in \Lambda$ and all
$t\in J$ that the dynamic ray $g_{\sigma^N\s}^{\kappa}$ is defined
at the potential $F^{\circ N}(t)$ with
\begin{eqnarray}
\Re\left(g_{\sigma^N\s}^{\kappa}(F^{\circ N}(t))\right)&\ge&
F^{\circ N}(t)-K-\left|r_{\kappa,\sigma^N\s}(F^{\circ N}(t))\right|>\n\\
&>&F^{\circ N}(t)-K-5 > 0\;.\label{notzero}
\end{eqnarray}
This shows the first statement.

For the second claim, let $(t_0,\kappa_0)$ be a pair of a
potential $t_0>t_{\s}$ and a parameter $\kappa_0\in\C$ such that
$\kappa_0\in \D_{\s}(t_0)$. Suppose by way of contradiction that
there are sequences $(t_n)_{n\ge 1}$ and $(\kappa_n)_{n\ge 1}$
with $t_n\to t_0$ and $\kappa_n \to\kappa_0$, such that
\[
\forall n\ge 1\ \ \exists N_n\ge 1\;,\ \exists t_n'\ge F^{\circ
N_n}(t_n):\quad g_{\sigma^{N_n}\s}^{\kappa_n}(t_n')=0\;.
\]
By the first part above, we may pass to a subsequence so that all
the $N_n$ are equal to some $N$. Furthermore, the sequence
$(t_n')_{n\ge 1}$ is contained in some compact interval $[F^{\circ
N}(t_*),t^*]$, where $t_*=\inf_n t_n'$ and $t^*$ is some potential
where we have good control, compare (\ref{notzero}). So by passing
to a subsequence once more we may assume that $(t_n')_n$ converges
to some $t_0'\ge F^{\circ N}(t_0)$. Note that
$\kappa_0\in\D_{\s}(t_0)$ implies $\kappa_0\in
\D_{\sigma^N\s}(t_0')$. So since the map $(t,\kappa)\mapsto
g_{\sigma^{N}\s}^{\kappa}(t)$ is continuous wherever it is
defined, it follows from $g_{\sigma^{N}\s}^{\kappa_n}(t_n')=0$ for
all $n\ge 1$ that $\lim_{n\to\infty}
g_{\sigma^{N}\s}^{\kappa_n}(t_n')= 0$, and therefore
$g_{\sigma^{N}\s}^{\kappa_0}(t_0')= 0$. This contradicts the
assumption $\kappa_0\in \D_{\s}(t_0)$. \qed
\end{proof}

\begin{proposition}[Discreteness and Local Cont.~Extension
of $G_{\s}$]\label{discreteness}
\lineclear Consider a
sequence $\s\in\S_0$ and a potential $t>t_{\s}$.
\begin{enumerate}
\item If $t_n\to t$ and $\kappa_n\to\kappa$ are sequences such
that $\kappa_n\in\G_{\s}(t_n)$ for all $n\ge 1$, then
$\kappa\in\G_{\s}(t)$. In particular, $\G_{\s}(t)$ is closed.

\item The set $\G_{\s}(t)$ is discrete in $\C$.

\item For every $\kappa_0\in \G_{\s}(t)$ there are neighborhoods
$\Lambda\subset\C$ and $J \subset\R$ containing $\kappa_0$ and $t$
respectively, such that for every $t'\in J$, the number of
elements of $\G_{\s}(t')\cap \Lambda$ (counting multiplicities)
equals the finite multiplicity of $\kappa_0$ as a root of the map
$\kappa\mapsto g_{\s}^{\kappa}(t)$.

More precisely, for every sequence $t_n\to t$ there is an $N\in\N$
and a sequence $(\kappa_n)_{n\ge N}\to\kappa_0$ such that
$\kappa_n\in\G_{\s}(t_{n})$ for all $n\ge N$.

\end{enumerate}
\end{proposition}

\begin{proof}
Let $t_n\to t$ and $\kappa_n\to\kappa$ be sequences such that
$\kappa_n\in\G_{\s}(t_n)$ for all $n\ge 1$. We have to show that
$\kappa\in\G_{\s}(t)$. By Lemma \ref{defdomain} (\ref{existsN}) we
may assume without loss of generality that there is an $N\in\N$
such that for all $m\ge N$ we have
$\{\kappa,\kappa_1,\kappa_2,\ldots\}\subset
\D_{\sigma^m\s}(F^{\circ m}(t_n))$. For every $n$ we have
$g_{\sigma^N\s}^{\kappa_n}(F^{\circ N}(t_n))=E_{\kappa_n}^{\circ
N}(0)$ and thus by continuity $g_{\sigma^N\s}^{\kappa}(F^{\circ
N}(t))=E_{\kappa}^{\circ N}(0)$. If $\kappa\not\in \D_{\s}(t)$,
then there is an $m\ge 1$ and a $t_0\ge F^{\circ m}(t)$ such that
$g_{\sigma^m\s}^{\kappa}(t_0)=0$ and thus
\[
E_{\kappa}^{\circ N}(0)=g_{\sigma^{m+N}\s}^{\kappa}(F^{\circ
N}(t_0))\;.
\]
We get two different potentials for $E_{\kappa}^{\circ N}(0)$,
which contradicts injectivity of $g^{\kappa}$ in Theorem
\ref{dynrays}. Therefore $\kappa_0\in\D_{\s}(t)$, and by
continuity $\kappa_0\in\G_{\s}(t)$.

For discreteness, consider a parameter $\kappa_0\in\G_{\s}(t)$ and
suppose that $\kappa_0$ is not isolated in $\G_{\s}(t)$. Let $U$
be the connected component of $\D_{\s}(t)$ which contains
$\kappa_0$. Since $\D_{\s}(t)$ is open, $U$ is open in $\C$.
Analyticity of $\kappa\mapsto g_{\s}^{\kappa}(t)$ and the identity
principle imply $U\subset \G_{\s}(t)\subset \D_{\s}(t)$. Therefore
$U$ is also the connected component of the closed set $\G_{\s}(t)$
containing $\kappa_0$ and therefore closed in $\C$. We conclude
$U=\C$, which is a contradiction.

For the third claim, consider neighborhoods $J_0,\Lambda_0$ of
$t,\kappa_0$ respectively as provided by Lemma \ref{defdomain}(2).
Since $\G_{\s}(t)$ is discrete, there is an $\eps>0$ such that
$\kappa_0$ is the only root of the map $\kappa\mapsto
g_{\s}^{\kappa}(t)$ within $\Lambda:=D_{\eps}(\kappa_0)$ and such
that $\ovl\Lambda\subset \Lambda_0$. Let $\gamma(s):=\kappa_0+\eps
e^{is}$. By Rouch\'e's Theorem, the multiplicity of $\kappa_0$ as
a zero equals the winding number $\eta(g_{\s}^{\gamma}(t),0)$ of
$g_{\s}^{\gamma}(t)$ around $0$. The holomorphic family
$\{f_{t'}(\kappa):=g_{\s}^{\kappa}(t'): \kappa\in \Lambda, t'\in
J_0\}$ is bounded and thus normal by Montel's Theorem. Hence if
$t_n\to t$ then $f_{t_n}$ converges uniformly to $f_t$ on
$\ovl\Lambda$. In particular we have
$\eta(g_{\s}^{\gamma}(t'),0)=\eta(g_{\s}^{\gamma}(t),0)$ for
potentials $t'$ sufficiently close to $t$.

By uniform convergence we can shrink $\gamma$ as $t_n$ gets closer
to $t$, and we find such parameters $\kappa_n$ which converge to
$\kappa_0$ as claimed in the additional statement. \qed
\end{proof}

We are now ready to state and to prove the main result.
\begin{theorem}[Parameter Rays at Their Full Length]\label{fulllength}
\lineclear For every external address $\s\in\S_0$ there is a
unique curve $G_{\s}:(t_{\s},\infty)\to I\subset\C$, called
\emph{parameter ray}, such that for all $t>t_{\s}$,
$\kappa_0=G_{\s}(t)$ satisfies $g_{\s}^{\kappa_0}(t)=0$ and is a
simple root of the map $\kappa\mapsto g_{\s}^{\kappa}(t)$.
Conversely, if $g_{\s}^{\kappa}(t)=0$ with $t>t_{\s}$, then
$\kappa=G_{\s}(t)$. The parameter rays are injective and pairwise
disjoint, and they carry the asymptotics
\[
G_{\s}(t)=t+ 2\pi i s_1 +R_{\s}(t) \quad \mbox{with} \quad
|R_{\s}(t)| =O(te^{-t}) \mbox{ as }t\to\infty\ .
\]
\end{theorem}

\begin{proof}
By Proposition \ref{parabound} and Lemma \ref{discreteness}(2),
$\G_{\s}(t)$ is bounded and discrete, and thus finite. Set
$n(t):=|\G_{\s}(t)|$ (counting multiplicities) and let
\[
J_{\s}:=\{T>t_{\s}:\ n(t)\ge 1 \ \forall t\ge T\}\;.
\]
By Corollary \ref{parendunique} we have $n(t)=1$ for every
sufficiently large $t$. The set $J_{\s}$ is thus non-empty, and it
follows from Proposition \ref{discreteness}(3) that $J_{\s}$ is
open. We will now show that $J_{\s}$ is also closed in
$(t_{\s},\infty)$. Let $t_*:=\inf J_{\s}$ and suppose
$t_*>t_{\s}$.  A function $t\mapsto G_{\s}(t)\in\G_{\s}(t)$ can be
defined on $J_{\s}$, possibly involving a choice. By Proposition
\ref{parabound}, the set $\{G_{\s}(t):\ t\in(t_*,t_*+1)\}$ is
contained in a compact set. Thus the set $L$ of all limits
$\lim_{t\searrow t_*} G_{\s}(t)$ is a nonempty compact subset of
$\C$. By Proposition \ref{discreteness}(1), it follows that
$L\subset \G_{\s}(t)$. Hence $J_{\s}$ is closed in
$(t_{\s},\infty)$ and $J_{\s}=(t_{\s},\infty)$.

Similarly we show that the set $J_{\s}':=\{t>t_{\s}:\ n(t)\ge 2\}$
is vacuous: by using Proposition \ref{discreteness} like above,
this set is open and closed relative $(t_{\s},\infty)$. However,
the complement is non-empty, because it contains an interval of
the form $(T,\infty)$ (on which Corollary \ref{parendunique}
holds). Therefore $J_{\s}'=\emptyset$, and for every $t>t_{\s}$ we
have $n(t)=1$. This shows that the choice for
$G_{\s}:(t_{\s},\infty)\to I$ above was unique.

Now this means that $G_{\s}$ is continuous because of Proposition
\ref{discreteness}(3). Injectivity of $G_{\s}$ and disjointness of
the parameter rays follow from the injectivity of $g^{\kappa}$ in
Theorem \ref{dynrays}, since every parameter $\kappa$ has at most
one external address and one potential. The asymptotic behavior
follows from Proposition \ref{parrayend}. \qed
\end{proof}

\remark Note that unlike dynamic rays, the parameter rays are
always defined on the entire interval $(t_{\s},\infty)$.

In \cite{frs}, the above result will be extended to endpoints for
fast sequences $\s$. This will yield a complete classification of
escaping parameters: there is a continuous bijection $G:X\to I$,
where $X$ is defined as in Section \ref{sec_dynrays}, and the
path-connected components of $I$ are exactly the parameter rays,
including the endpoints at fast addresses. Moreover, one can
easily show \cite{f} that the parameter rays are $C^1$-curves, and
it seems that with some more work one can also show $C^{\infty}$.

\subsection{Vertical Order of Parameter Rays}

We show that parameter rays have a natural vertical order which
coincides with the lexicographic order of their external addresses.

\begin{definition}[Vertical Order]
\lineclear Let $\Gamma$ be a family of injective rays
$\gamma:\R^+\to\C$ with $\Re\gamma(t) \to \infty$ as $t\to\infty$,
such that any two curves $\gamma_1\neq\gamma_2\in\Gamma$ have a
bounded set of intersection. For $\gamma_1,\gamma_2\in\Gamma$ let
$\H_R$ be a right half-plane in which $\gamma_1$ and $\gamma_2$
are disjoint, and denote $\H_R^+(\gamma_2)$ and $\H_R^-(\gamma_2)$
the upper respectively lower component of $\H_R\setminus
\gamma_2$. Then
\[
\gamma_1\succ\gamma_2:\quad\Longleftrightarrow\quad
\gamma_1\subset H_R^+(\gamma_2)
\]
defines a linear order on $\Gamma$. We say that $\gamma_1$ is
\emph{above} $\gamma_2$. \qed
\end{definition}

Equip $\S=\Z^{\N^*}$ with the lexicographic order `$>$'.

\begin{lemma}[Vertical Order of Dynamic
Rays]\label{vert_order_dyn}
\lineclear For all $\kappa\in\C$ and all exponentially bounded
addresses $\s,\tilde\s\in\S_0$,
\[
g_{\s}^{\kappa}\succ g_{\tilde \s}^{\kappa}
\quad\Longleftrightarrow \quad \s > \tilde\s\;.
\]
\end{lemma}

\begin{proof}
Without loss of generality, say $\s>\tilde \s$. If $s_1>\tilde
s_1$, then $g_{\s}^{\kappa}\succ g_{\tilde\s}^{\kappa}$ follows
directly from the asymptotic estimate (\ref{rK_est}) in Theorem
\ref{dynrays}. Otherwise let $k> 1$ be the first entry in which
$\s$ and $\tilde \s$ differ. Then by the same argument,
$g_{\sigma^{k-1}\s}^{\kappa}\succ
g_{\sigma^{k-1}\tilde\s}^{\kappa}$. Since
$$g_{\s}^{\kappa}(t)=L_{\kappa,s_1}\circ \ldots \circ
L_{\kappa,s_{k-1}}\circ
g_{\kappa,\sigma^{k-1}\s}^{\kappa}(F^{\circ (k-1)}(t))\;,$$ where
the translated logarithms $L_{\kappa,s}=\Log-\kappa +2\pi i s$
preserve the vertical order, the claim follows. \qed
\end{proof}

\begin{proposition}[Vertical Order of Parameter Rays]
\lineclear For all exponentially bounded addresses
$\s,\tilde\s\in\S_0$,
\[
G_{\s}\succ G_{\tilde \s} \quad\Longleftrightarrow \quad \s >
\tilde\s\;.
\]
\end{proposition}

\begin{proof}
The claim follows from the asymptotic estimate in Theorem~\ref{fulllength} if the
first entries in $\s$ and $\tilde \s$ are different, so we assume they are equal
and $\s>\tilde \s$.

Consider the external address $\s=s_1s_2s_3\dots$ and set
\[
S:=\{z\in\C\colon \Re(z)\ge \xi, 2\pi (s_1-2)\le \Im(z)\le 2\pi (s_1+2)\}
\]
depending on $\xi\in\R$. By Theorem~\ref{dynrays}, part 2~(b), we can fix $\xi$ so
that for all $\kappa\in S$, there is a potential $\tau_\kappa>t_\s$ such that 
$g_{\s}^{\kappa}(t)$ is defined for $t\geq \tau_\kappa$ and 
$\Re(g_{\s}^{\kappa}(\tau_\kappa))<-11$.

In the dynamical plane for $\kappa\in S$, consider the right half plane
\[
H_\kappa:=\{z\in\C\colon \Re(z)>\Re(g_{\s}^{\kappa}(\tau_\kappa))\}
\,\,.
\]
The ray tail $g_{\s}^{\kappa}((\tau_\kappa,\infty))$ cuts $H_\kappa$ into two
unbounded components $H^+_\kappa$ and $H^-_\kappa$  (plus possibly some bounded
components) so that in $H^+_\kappa$, imaginary parts are unbounded above, while in
$H^-_\kappa$, they are unbounded below. The asymptotics in (\ref{rK_est}) from
Theorem~\ref{dynrays} together with the condition
$\Re(g_{\s}^{\kappa}(\tau_\kappa))<-11$
implies that for every $\kappa\in S$, either $0\in H^+_\kappa\cup H^-_\kappa$ or
$0=g_{\s}^{\kappa}(t)$ for some $t\geq\tau_\kappa$.

Define $S^+_0:=\{\kappa\in S\colon 0\in H^+_\kappa\}$ and analogously $S^-_0$. By
construction, $S\subset S^+_0\dot\cup S^-_0\dot\cup G_\s((t_\s,\infty))$.

The asymptotics of the parameter ray $G_\s$ from Theorem~\ref{fulllength} implies
that $S\sm G_\s((t_s,\infty))$ contains two unbounded components (plus possibly
some bounded ones); let $S^+$ and $S^-$ be the unbounded component above resp.\
below $G_\s((t_s,\infty))$ in the obvious sense.

We claim that $S^+\subset S^+_0$ and $S^-\subset S^-_0$. Indeed, the set $S^+$
contains a tail of the parameter ray $G_{\s'}$ with $\s'=(s_1+1)s_2s_3s_3\dots$.
For parameters $\kappa$ on this tail, the vertical order of dynamic rays implies
$0\in H^+_\kappa$, hence $\kappa\in S^+_0$. Since $S^+$ is connected, it follows
that $S^+\subset S^+_0$, and analogously $S^-\subset S^-_0$. 

The parameter ray $G_{\tilde \s}$ also has a tail in $S$, and for parameters
$\kappa$ on this tail, the vertical order of dynamic rays implies
$\kappa\in S^-_0$. If $\Re(\kappa)$ is sufficiently large, then $\kappa\in S^+\cup
S^-$. Finally, $\kappa\in S^+$ would imply $\kappa\in S^+_0$, a contradiction, so
we conclude $\kappa\in S^-$.
\qed
\end{proof}

\hide{
\begin{proof}
If $s_1\neq \tilde s_1$, the claim follows directly from the
asymptotic estimate in Theorem \ref{fulllength}. So assume that
$s_1=\tilde s_1$ and $\s > \tilde \s$. Consider the function
$f=g^{\bullet}_{\s}(t):\kappa\mapsto g_{\s}^{\kappa}(t)$. Let us
only consider sufficiently large potentials $t>T$, where the
parameter rays $G_{\s}$, $G_{\tilde\s}$ are close to straight
lines, such that $f$ is defined for all $\kappa\in
\Lambda:=B_{3\pi}(t+2\pi is_1)$, i.e. $f:\Lambda\to\C$ is a
continuous function in $\kappa$ (and $t$). In fact,
$f(\kappa)=(t+2\pi is_1) -\kappa+r_{\kappa,\s}(t)$ is a linear
function up to an error of $|r_{\kappa,\s}(t)|=O(|\kappa|e^{-t})$.
Since $t$ is large, we may assume that $G_{\tilde
\s}(t),G_{\tilde\s}(t)+2\pi i \in \Lambda$.
Let $\eps =2R$, where $R$ is an estimate for $|r_{\kappa,\s}(t)|$
on $\Lambda$. Consider the `square' $Q\subset \Lambda$ with
horizontal edges $G_{\tilde\s}(t+t')$ and $G_{\tilde\s}(t+t')+2\pi
i$ ($t'\in [-\eps,\eps]$) and straight vertical edges. Lemma
\ref{vert_order_dyn} yields that $\Im f <0$ on the lower edge and
$\Im f>0$ on the upper edge. Moreover, $Re f$ is negative
(positive) on the right (left) edge. By continuity, the only zero
$G_{\s}(t)$ of $f$ is contained in $Q$ and thus above
$G_{\tilde\s}(t)$. Since the parameter rays are parameterized by
real parts (up to a negligible error), the claim follows.
\qed
\end{proof}
}

\section{The Proof of the Bound on Parameter
Rays}\label{sec_bound}

\subsection{First Derivative of Dynamic Rays}\label{secderiv}
In order to prove Proposition \ref{parabound}, we will need
estimates on the derivative of dynamic rays, which lead to
estimates on the winding numbers of dynamic rays. These in turn
will help us control the rays at small potentials in order to
obtain a bound on the absolute value of all $\kappa\in \G_{\s}(t)$
with prescribed combinatorics $(\s,t)$.

We will often be concerned with obtaining estimates on some ``tail
pieces'' of dynamic rays. In order to simplify the discussion
without having to keep track of exact constants, we make the
following definition.

\begin{definition}[Properties on Ray Tails]
\label{raytails}  \lineclear 
We say that a property $P(\kappa,\s,t)$ holds \emph{on ray tails} if there are
$A,B,C\ge 0$ such that for all $\s\in\S_0$ and all $K\ge 1$, the property
$P(\kappa,\s,t)$ holds whenever $|\kappa| \le K$ and
\[
t\ge At^*_\s+B\log K+C
\;,
\]
where $t^*_\s$ is the constant from Definition and Lemma
\ref{minpot}. 
\hide{
$t\ge T(K,\s)$,
where $T:[1,\infty)\times \S_0\to \R^+$ is defined by
\[
T(K,\s):=At_{\s}^*+B\log K+C 
\]
and $t_{\s}^*$ is the constant from Definition and Lemma
\ref{minpot}. Moreover, let
\[
n_{\s}^T(K,t):=\min\{n\ge 0:\ F^{\circ n}(t)\ge
T(K,\sigma^n\s)\}\;.
\]
}
\end{definition}

\remark 
The problem is that it is much easier to control tails of dynamic
rays than to control entire rays. This control is non-uniform in $\kappa$: if
$|\kappa|$ is large, then we have good control only for large potentials $t$. The
following result often allows to transfer results from ray tails to all rays.

\begin{lemma}[From Ray Tails to Entire Rays]
\label{Lem:FromTailsToRays} \lineclear
Suppose a property $P(\kappa,\s,t)$ holds on ray tails and is backward invariant,
i.e.\ it holds for $g_{\kappa,\s}(t)$ whenever it holds for
$g_{\kappa,\sigma(\s)}(F(t))$. Then it holds on all dynamic rays.
\end{lemma}
\proof
The property holds on $g_{\kappa,\s}(t)$ as soon as there is an $n\in\N$ 
such that it holds on $g_{\kappa,\sigma^n(\s)}(F^{\circ n}(t))$. This is true as
there is an $N$ such that $F^{\circ N}(t)\ge At^*_{\sigma^N(\s)}+B\log K+C$ and
this follows from the last claim in Definition and Lemma \ref{minpot}.
\qed

The quantifier ``on ray tails'' commutes with finite,
but not infinite conjunctions: ``$\forall n:\, P_n(\kappa,\s,t)$
on ray tails'' is weaker than ``on ray tails, $\forall n:\,
P_n(\kappa,\s,t)$'': in the first case, the constants $A$, $B$, $C$ may depend on
$n$.

Using this notation, we can now say that the asymptotic bound
(\ref{rK_est}) of Theorem \ref{dynrays} holds on ray tails: in this case, for all
$t\ge t_\s^*+2\log K+2\log 3 \ge t_\s^*+2\log(K+3)=t_\s^K$.
This gives us very good control on the orbit
of points on dynamic rays, except for at most finitely many steps.
The following lemma helps to estimate after how many iteration steps good control
takes over.

\begin{lemma}[Bound on Initial Iteration Steps]
\label{Lem:BoundInitialSteps} \lineclear
Fix $\s\in\S_0$, $t>t_{\s}$, $A\geq 1$ and $B,C\geq 0$. Then for every $K\ge 1$,
there is an $N\in\N$ such that for all $n\geq N$, we have
\begin{equation}
F^{\circ n}(t)\ge At_{\sigma^n(\s)}^*+B\log K+C \;.
\label{Eq:Bound_on_N}
\end{equation}
The value of $N$ has the following properties for fixed $\s$:
\begin{enumerate}
\item
for fixed $K$, it is (weakly) monotonically decreasing in $t$;
\item
if $t\ge At_{\s}^*+B\log K+C$, then $N=0$;
\item
if for fixed $t$, $N_0$ is such that
\[
F^{\circ N_0}(t)\ge At_{\sigma^{N_0}(\s)}^*+C+1 \;,
\]
then (\ref{Eq:Bound_on_N}) holds for $N=N_0+N_1$ as soon as $F^{\circ N_1}(1)\ge
B\log K$.
\end{enumerate}
\end{lemma}
\proof
Note first that by convexity, $F(x+y)\ge F(x)+F(y)$ for all $x,y\ge 0$. Similarly,
$F(Ax)\ge AF(x)$ and in particular $F(x)\ge x$. Moreover, by definition,
$F(t^*_{\s})\ge t^*_{\sigma(\s)}$. This implies that if (\ref{Eq:Bound_on_N}) holds
for $n$, then it also holds for $n+1$. The second claim follows, and the first is
trivial.

The third claim is verified as follows:
\begin{eqnarray*}
F^{\circ N}(t)=F^{\circ(N_0+N_1)}(t) &\ge&
F^{\circ N_1}(At_{\sigma^{N_0}(\s)}^*+C+1) 
\\
&\ge&
F^{\circ N_1}(At^*_{\sigma^{N_0}(\s)})+   F^{\circ N_1}(C)+   F^{\circ N_1}(1)
\\
&\ge&
At^*_{\sigma^N(\s)} + C + B\log K \;.
\end{eqnarray*}
\qed

Using the terminology of ``properties on ray tails'', the following statements
follow easily from \cite{sz}, Lemma $3.3$ and Proposition $3.4$.

\begin{lemma}[Further Properties of Dynamic Rays]\label{propdynrays}
\lineclear Let $\s\in\S_0$, $K>0$ and $\kappa$ be a parameter with
$|\kappa|\le K$. On the interval $(t_{\s}^K,\infty)$, the curve
$g_{\s}^{\kappa}$ is the uniform limit of the functions
$g_{\kappa,\s}^n$ defined by $g_{\kappa,\s}^0:=\mathrm{id}$ and
$g_{\kappa,\s}^{n+1}(t):=L_{\kappa,s_1}\circ
g_{\kappa,\sigma\s}^n(F(t))$, i.e.
\begin{equation}\label{gn}
g_{\kappa,\s}^n(t):=L_{\kappa,s_1}\circ\cdots\circ
L_{\kappa,s_n}(F^{\circ n}(t))\;.
\end{equation}
On ray tails, they satisfy for all $n,k\ge 1$
\begin{eqnarray}
|g_{\kappa,\s}^n(F^{\circ k}(t))|&\ge& 2^k\quad\mbox{and}\\
|g_{\kappa,\s}^k(t)-g_{\s}^{\kappa}(t)|&\le & 2^{-k}\;.
\end{eqnarray}
\qed
\end{lemma}

We omit the straightforward but technical proof of the following
lemma.

\begin{lemma}[Some Properties of $F$]\label{PropF}
\lineclear If $x\ge 0$ and $t\ge 5$ are real numbers such that
$t\ge 2x+5$ then
\[
\sum_{k=1}^{\infty}\frac{1}{F^{\circ
k}(t)+1}<3e^{-t}\quad\mbox{and}\quad
\sum_{k=1}^{\infty}\frac{F^{\circ k}(x)}{F^{\circ k}(t)+1} <
(e^x+1)e^{-t}\;.
\]
Moreover,
\begin{eqnarray}
\frac{d}{dt}F^{\circ n}(t)&=& \prod_{k=1}^n (F^{\circ
k}(t)+1)\quad\mbox{and}\label{ddtFn}\\
\exists T>0:\ \forall t\ge T,n\ge 1:\quad && \frac{(F^{\circ
n})'(t)}{(F^{\circ n}(t)+1)^2}\le \frac1{F^{\circ
n}(t-1)+1}\;.\label{Fk'Fk2}
\end{eqnarray}
%
\qed
\end{lemma}

Differentiability of dynamic rays has already been proven in 1988
by M.~Viana da Silva \cite{v}. We will prove it again in order to
obtain explicit estimates on the first and second derivatives.

\begin{proposition}[The Derivative of Dynamic Rays]\label{ddtdynray}
\lineclear For every sequence $\s\in\S_0$ and every parameter
$\kappa$, the dynamic ray $g_{\s}^{\kappa}(t)$ is continuously
differentiable with respect to the potential $t$ with derivative
\begin{equation}\label{ddteqn}
(g_{\s}^{\kappa})'(t)=\prod_{k=1}^{\infty} \frac {F^{\circ
k}(t)+1}{g_{\sigma^k\s}^{\kappa} (F^{\circ k}(t))} \,\neq 0\,\,\,.
\end{equation}
Moreover, on ray tails,
\begin{equation}\label{gprimeest}
|(g_{\s}^{\kappa})'(t)-1|< e^{-t/2} \,\,.
\end{equation}
\end{proposition}

\begin{proof}
It is sufficient to prove that $g_{\s}^{\kappa}$ is differentiable
and satisfies (\ref{ddteqn}) and (\ref{gprimeest}) on ray tails:
if (\ref{ddteqn}) is known for $g_{\sigma\s}^{\kappa}$ at $F(t)$,
then
\begin{eqnarray*}
(g_{\s}^{\kappa})'(t)&=&\frac{(E_{\kappa}\circ
g_{\s}^{\kappa})'(t)}{E_{\kappa}(g^{\kappa}_{\s}(t))}=\frac{(g_{\sigma\s}^{\kappa}\circ
F)'(t)}{g^{\kappa}_{\sigma\s}(F(t))}=\\
&=&\left(\prod_{k=1}^{\infty}\frac{F^{\circ
(k+1)}(t)+1}{g^{\kappa}_{\sigma^{k+1}}(F^{\circ
(k+1)}(t))}\right)\frac{F(t)+1}{g^{\kappa}_{\sigma\s}(F(t))}=\prod_{k=1}^{\infty}
\frac {F^{\circ k}(t)+1}{g_{\sigma^k\s}^{\kappa} (F^{\circ
k}(t))}\;.
\end{eqnarray*}
By Lemma~\ref{Lem:FromTailsToRays}, (\ref{ddteqn}) holds on all rays.

Recall the functions $g_{\kappa,\s}^n(t)$, defined in Lemma
\ref{propdynrays}, which converge uniformly to
$g_{\s}^{\kappa}(t)$. By the chain rule, for every $\kappa\in\C$,
$\s\in\S_0$, $n\ge 1$ and $t$ (where defined),
\begin{eqnarray*}
\frac{d}{dt}L_{\kappa,s_1}\circ\cdots\circ
L_{\kappa,s_n}(t)&=&\prod_{k=1}^{n}\left(g_{\kappa,\sigma^k\s}^{n-k}(F^{\circ
(k-n)}(t))\right)^{-1}\;,
\end{eqnarray*}
and thus together with (\ref{ddtFn}) in Lemma \ref{PropF}, by the
chain rule again,
\begin{eqnarray}
(g^n_{\kappa,\s})'(t)=\left(\prod_{k=1}^{n} \frac {F^{\circ
k}(t)+1} {g_{\sigma^k\s}^{\kappa}(F^{\circ k}(t))}\right) \cdot
\left(\prod_{k=1}^{n}
\underbrace{\frac{g_{\sigma^k\s}^{\kappa}(F^{\circ k}(t))}
 {g^{n-k}_{\kappa,\sigma^k\s}(F^{\circ k}(t))}}_{=:P_k^n(t)}\right)
\,\,.\label{gnprime}
\end{eqnarray}
Let us first show that (on ray tails) $\prod_{k=1}^n P_k^n(t)$
converges uniformly to $1$ as $n\to\infty$. Indeed, on ray tails,
\[
|P_k^n(t)-1|=\left|\frac{g_{\sigma^k\s}^{\kappa}(F^{\circ
k}(t))-g^{n-k}_{\kappa,\sigma^k\s}(F^{\circ
k}(t))}{g^{n-k}_{\kappa,\sigma^k\s}(F^{\circ k}(t))} \right|\le
\frac 1{2^{n-k}}\cdot\frac1{2^k}=\frac1{2^n}
\]
by Lemma \ref{propdynrays}. Thus $1-(1-2^{-n})^n\le |\prod_{k=1}^n
P_k^n(t)-1|\le (1+2^{-n})^n-1$, which means that $\prod_{k=1}^n
P_k^n(t)$ converges to $1$ uniformly in $t$.

By Weierstra{\ss}' Theorem, it only remains to show that the first
product of (\ref{gnprime}) converges uniformly on ray tails and
satisfies the uniform bound (\ref{gprimeest}) there. Note that
\[
\Log \prod_{k=1}^{n} \frac{F^{\circ
k}(t)+1}{g_{\sigma^k\s}^{\kappa} (F^{\circ k}(t))}=-\sum_{k=1}^{n}
\Log\left(1+\underbrace{\frac{g_{\sigma^k\s}^{\kappa} (F^{\circ
k}(t))-F^{\circ k}(t)-1}{F^{\circ k}(t)+1}}_{=:x_k(t)}\right)\;.
\]
Again by Lemma \ref{propdynrays}, on ray tails,
$$|x_k(t)|\le
\frac{|\kappa|+2\pi |s_{k+1}|+2}{F^{\circ k}(t)+1}\le
\frac{|\kappa|+2\pi F^{\circ k}(t_{\s}^*)+2}{F^{\circ k}(t)+1}\le
1/2\;.
$$
Thus on ray tails, by the first inequality from Lemma \ref{PropF},
\[
\sum_{k=1}^{\infty}|x_k(t)|<\left(3(|\kappa|+2)+2\pi(e^{t_{\s}^*}+1)
\right)e^{-t}\le e^{-t/2}/4\;.
\]
Since $|\Log(1+x)|\le 2|x|$ for $|x|\le 1/2$, it follows on ray
tails:
\begin{eqnarray*}
\left|\Log \prod_{k=1}^{\infty} \frac{F^{\circ
k}(t)+1}{g_{\sigma^k\s}^{\kappa} (F^{\circ k}(t))}\right|\le
\sum_{k=1}^{\infty}|\Log(1+x_k(t))|
\le\sum_{k=1}^{\infty}2|x_k(t)|\le e^{-t/2}/2\le 1/2\;.
\end{eqnarray*}
Finally, (\ref{gprimeest}) follows on ray tails, since $|z-1|\le
2|\log z|$ for $|\log z|\le 1/2$.
 \qed
\end{proof}

\subsection{Second Derivative of Dynamic Rays}

\begin {proposition} [The Second Derivative of Dynamic Rays]
\label{d2dtdynray}\lineclear 
Every dynamic ray $g^{\kappa}_{\s}:(t_{\s},\infty)\to
\C$ is twice continuously differentiable. On ray tails,
\begin{equation}
\label{est_2nd_der}  |(g_{\s}^{\kappa})''(t)| < e^{-t/2}
\quad\text{and}\quad \left |
 \frac {(g^{\kappa}_{\s})''(t)}{(g^{\kappa}_{\s})'(t)}\right| < e^{-t/2}\ .
\end{equation}
\end {proposition}

\begin{proof}
Define $f_{\s}(t):=(t+1)/g_{\s}^{\kappa}(t)$ and for $k\ge 1$
\[
P_k(t):=\frac {F^{\circ k}(t)+1}{g_{\sigma^k\s}^{\kappa} (F^{\circ
k}(t))}=f_{\sigma^k\s}(F^{\circ k}(t))\;. \] The partial products
$h_N(t):=\prod_{k=1}^N P_k(t)$ in (\ref{ddteqn}) from Proposition
\ref{ddtdynray} converge uniformly to $(g_{\s}^{\kappa})'$ and we
thus need to show that the derivatives $h_N'$ converge uniformly
and that the limit satisfies (\ref{est_2nd_der}). By the product
rule,
\[
h_N'(t)=h_N(t)\cdot\sum_{k=1}^N\frac{P_k'(t)}{P_k(t)}\;.
\]
The factor $h_N(t)$ can be bounded by $2$. (This also shows that
the first estimate in (\ref{est_2nd_der}) implies the second one.)
Since $|s_{k+1}|\le F^{\circ k}(t_{\s}^*)$, on ray tails,
\[
|P_k(t)|\ge \frac{F^{\circ k}(t)+1}{F^{\circ k}(t)+|\kappa|+2\pi
|s_{k+1}|+1}\ge 1/2\quad \mbox{for every }k\ge 1\;.
\]
It remains to show that $\sum P_k'(t)$ converges on ray tails.
We estimate
\begin{eqnarray*}
|f_{\sigma^k\s}'(F^{\circ
k}(t))|&=&\left|\frac{g_{\sigma^k\s}^{\kappa}(F^{\circ
k}(t))-(F^{\circ k}(t)+1)
(g_{\sigma^k\s}^{\kappa})'(F^{\circ k}(t))}{(g_{\sigma^k\s}^{\kappa}(F^{\circ k}(t)))^2}\right|=\\
&=& \left|\frac{-\kappa+2\pi i s_{k+1}-1+O(F^{\circ
k}(t)e^{-2F^{\circ k}(t)/3})}{(F^{\circ k}(t)-\kappa+2\pi
is_{k+1}+O(e^{-F^{\circ k}(t)}))^2}\right|\le\\
&\le& \frac{2\pi F^{\circ k}(t_{\s}^*)+|\kappa|+2}{((F^{\circ
k}(t)+1)/2)^2}
\end{eqnarray*}
and thus by (\ref{Fk'Fk2}) in Lemma \ref{PropF}
\[
|P_k'(t)|=|f_{\sigma^k\s}'(F^{\circ k}(t))|\cdot (F^{\circ
k})'(t)\le 4\cdot \frac{2\pi F^{\circ
k}(t_{\s}^*)+|\kappa|+2}{F^{\circ k}(t-1)+1}\;.
\]
Using Lemma \ref{PropF} once more,
\begin{eqnarray*}
\sum_{k=1}^{\infty}|P_k'(t)|&\le& \left(3(4|\kappa|+8)+8\pi
(e^{t_{\s}^*}+1)\right)e^{-t+1}\le e^{-t/2}/4
\end{eqnarray*}
on ray tails. This shows that the $h_N'$ converge uniformly, and
\[
|(g_{\s}^{\kappa})''(t)|\le |h_N'(t)|\le 2\sum_{k=1}^{\infty}
2|P_K'(t)|\le e^{-t/2}\;.
\]
\vspace{-1.2cm}

\qed
\end{proof}

\hide{
For easier reference, we rephrase part of the previous result as follows,
expanding Definition~\ref{raytails}.
\begin{corollary}[The Second Derivative of Dynamic Rays]
\label{Cor:SecondDerivative} \lineclear
For every $\s\in\S_0$, there is a function $(K,t_0)\mapsto n_{\s}(K,t_0)$
such that for all $|\kappa|\le K$, $n\ge n_{\s}(K,t_0)$ and $t\ge
t_0>t_{\s}$
\begin{equation}\label{Eq:d2dtest}
\left|\frac{g''_{\kappa,\sigma^{n}\s}(F^{\circ n}(t))}
{g'_{\kappa,\sigma^{n}\s}(F^{\circ n}(t))}\right|<e^{-F^{\circ
n}(t)/2} < 1\;.
\end{equation}
This function is monotonically increasing in $K$ and decreasing in
$t$
\hide{, so that for given $K$ and $\s$, we have $n_{\s}(K,t_0)=0$ for all
sufficiently large $t_0$
}
.
\qedd
\end{corollary}
}

\subsection{Variation Numbers of Dynamic Rays}\label{sec_winnr}

Several key ideas in this section are are due to Niklas Beisert.

For a closed $C^1$-curve $\gamma:[t_0,t_1]\to\C$ and
$a\in\C\sm\gamma([t_0,t_1])$, \emph{the winding number of $\gamma$
around $a$} is defined by
\[
\eta(\gamma,a):=\frac{1}{2\pi}\int_{t_0}^{t_1}
d\arg(\gamma(t)-a)\quad \in\Z\;.
\]

\begin{definition}[Variation
Number]\label{def_windingnr} \lineclear Consider a $C^1$-curve
$\gamma:(t_0,\infty)\to\C$ with $t_0\ge -\infty$ and $a\not\in
\gamma(t_0,\infty)$. Define the \emph{variation number of}
$\gamma$ \emph{around} $a$ by
\begin{eqnarray*}
\alpha(\gamma,a):=\frac{1}{2\pi}\int_{t_0}^\infty
\left|\Im\frac{\gamma'(t)}{\gamma(t)-a} \right|dt=
\frac{1}{2\pi}\int_{t_0}^\infty
\left|d\arg(\gamma(t)-a)\right|\quad\in [0,\infty]\;.
\end{eqnarray*}
\end{definition}
Unlike the winding number, the variation number also measures
local oscillations of the curve.

\begin{definition}[Admissible Curves]
\lineclear An \emph{admissible curve} is an injective $C^2$-curve
$\gamma:(t_0,\infty)\to\C$ with non-vanishing derivative $\gamma'$
such that for $t\to\infty$
\[
|\gamma(t)|=t+O(1)\;,\quad|\gamma'(t)|=1+O(1/t)\;,\quad\mbox{and}\quad
|\gamma''(t)|=O(1/t^2)\;.
\]
\end{definition}

\begin{lemma}[Admissible Curves Have Variation Numbers]
\lineclear Let $\gamma:(t_0,\infty)\to\C$ be an admissible curve.
\begin{enumerate}
\item If $a\in\C\sm\overline{\gamma(t_0,\infty)}$ and
$|\gamma'(t)|$ is bounded as $t\searrow t_0$, then
$\alpha(\gamma,a)$ is finite.

\item For every $t_1>t_0$,
$\alpha(\gamma|_{(t_1,\infty)},\gamma(t_1))$ is finite.

\item For every $t_1>t_0$, $\alpha(\gamma'|_{(t_1,\infty)},0)$ is
finite.
\end{enumerate}
\end{lemma}

\begin{proof}
In all statements the integrands are locally Riemann integrable.
Therefore we only have to show that the integrals
$\int_{t_0}^\infty |\Im\frac{\gamma'}{\gamma-a}|$ are finite near
the boundaries of integration.

Let us first discuss the lower boundary of integration. In Case 1,
we can bound $|\gamma(t)-a|$ below and $|\gamma'(t)|$ above. Fix
$\eps>0$. For Case 3 the continuous function
$|\gamma''(t)/\gamma'(t)|$ is bounded on $[t_1,t_1+\eps]$. In Case
2 however, the denominator tends to $0$. By the Taylor Theorem
applied to $\gamma\in C^2$, for every $t\in(t_1,t_1+\eps)$ there
is a $\xi\in [t_1,t]$ such that
\begin{eqnarray*}
\left|\Im \frac{\gamma'(t)}{\gamma(t)-\gamma(t_1)}\right|&=&
\left|\Im
\frac{\gamma'(t)}{\gamma'(t)(t-t_1)+\gamma''(\xi)(t-t_1)^2/2}\right|=\\
&=&\left|\Im\frac
1{t-t_1+\frac{\gamma''(\xi)}{2\gamma'(t)}(t-t_1)^2}\right|\;.
\end{eqnarray*}
For $\xi,t\in [t_1 , t_1+\eps]$,
$\frac{\gamma''(\xi)}{2\gamma'(t)}=:c(t)$ can be estimated
uniformly and is thus of class $O(1)$. Now for
$t-t_1=:\delta\searrow 0$ we observe
\[
\Im\frac 1{\delta+c\delta^2}=\frac{\delta^2 \Im \bar
c}{|\delta+c\delta^2|^2}=\frac{-\Re c}{|1+c\delta|^2}=O(1)\;.
\]
It is left to show that the limits $\lim_{x\to\infty}\int^x$ are
finite: for the first two cases we have (with $\gamma(t_1)=:a$)
\begin{eqnarray*}
\Im \frac{\gamma'(t)}{\gamma(t)-a}&=&
\Im\frac{1+O(t^{-1})}{t+O(1)}=\Im\frac{(1+O(t^{-1}))\overline{(t+O(1))}}{|t+O(1)|^2}=\\
&=&\frac{\Im(t+O(1))}{|t+O(1)|^2}=O(t^{-2})\;,
\end{eqnarray*}
and for the last case we estimate $|\Im
(\gamma''(t)/\gamma'(t))|\le |\gamma''(t)/\gamma'(t)|=O(1/t^2)$.
\qed
\end{proof}

\begin{lemma}[The Variation Number of a Half Line]\label{halfline}
\lineclear For the curve $\ell:\R^+\to\C$, $\ell(s)=\lambda s$
($\lambda\in\C^*$) we have for every $a\not\in \lambda\R^+_0$
\[
\alpha(\ell,a)=\frac{|\arg(a/\lambda)-\pi|}{2\pi}\;,
\]
if the argument is chosen in the interval $[0,2\pi]$.
\qed
\end{lemma}

The proof is left to the reader. The following lemma gives a very
useful connection between the variation number of a curve and of
its derivative.

\begin{lemma}[Variation Numbers of a Curve and its Derivative]\label{BoundVar}
\lineclear Let $\gamma_0:(t_0,\infty)\to\C$ be an admissible
curve, and $\gamma:=\gamma_0|_{(t_1,\infty)}$ its restriction for
some $t_1>t_0$. Then for every $a\not\in\gamma_0([t_1,\infty))$
\[
\alpha(\gamma,a)\le \alpha(\gamma',0)+1/2\ .
\]
\end{lemma}

\begin{figure}
\begin{picture}(0,0)%
\includegraphics{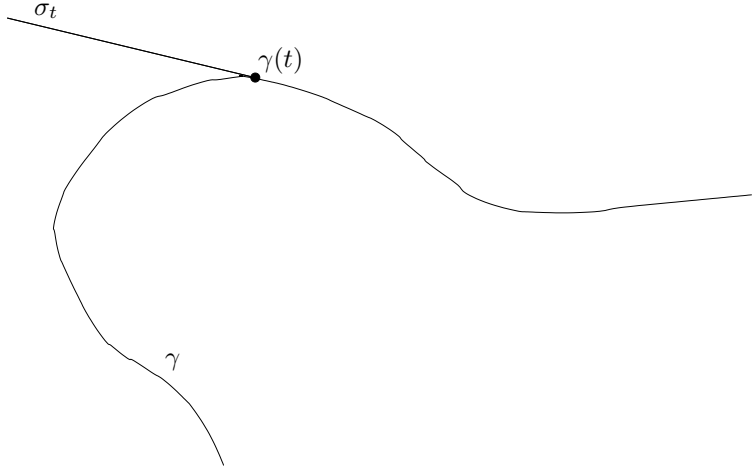}%
\end{picture}%
\setlength{\unitlength}{2486sp}%
\begingroup\makeatletter\ifx\SetFigFont\undefined%
\gdef\SetFigFont#1#2#3#4#5{%
  \reset@font\fontsize{#1}{#2pt}%
  \fontfamily{#3}\fontseries{#4}\fontshape{#5}%
  \selectfont}%
\fi\endgroup%
\begin{picture}(7447,4704)(1071,-5698)
\put(1351,-1186){\makebox(0,0)[lb]{\smash{\SetFigFont{10}{12.0}{\familydefault}{\mddefault}{\updefault}{\color[rgb]{0,0,0}$\sigma_t$}%
}}}
\put(3586,-1711){\makebox(0,0)[lb]{\smash{\SetFigFont{10}{12.0}{\familydefault}{\mddefault}{\updefault}{\color[rgb]{0,0,0}$\gamma(t)$}%
}}}
\put(2656,-4681){\makebox(0,0)[lb]{\smash{\SetFigFont{10}{12.0}{\familydefault}{\mddefault}{\updefault}{\color[rgb]{0,0,0}$\gamma$}%
}}}
\end{picture}
\caption{The linear continuation as defined in the proof of Lemma
\ref{BoundVar}.} \label{FigWinding}
\end{figure}

\begin{proof}
For every $t\ge t_1$ define the curve $\sigma_t:\R\to\C$ to be the
linear continuation of $\gamma$ at $t$, i.e.
$\sigma_t(s):=\gamma(s)$ if $s\ge t$ and
$\sigma_t(s):=\gamma(t)+(s-t)\gamma'(t)$ if $s\le t$, see Figure
\ref{FigWinding}. Consider the following two functions
$u,v:(t_1,\infty)\to\R$:
\begin{eqnarray*}
u(t)&=&\frac1{2\pi}\int^\infty_t |d\arg\gamma'| \qquad \mbox{and} \\
v(t)&=& \frac1{2\pi}\int^\infty_{-\infty} |d\arg(\sigma_t-a)|\;.
\end{eqnarray*}
Note that the form $d\arg$ is defined everywhere, although $\arg$
may not be.

Lemma \ref{halfline} (generalized to the degenerate case, which we
define via $\alpha(\ell,\ell(s)):=1/2$) gives us
\begin{eqnarray*}
\frac d{d t} \int_{-\infty}^t |d\arg(\sigma_t(s)-a)| &=& 2\pi
\frac d{dt} \alpha(\sigma_t|_{(-\infty,t)})=\frac
d{dt}\left|\arg\left(\frac{a-\gamma(t)}{\gamma'(t)}-\pi\right)\right|=\\
&=&\rho(t)\frac d{d
t}\left(\rule{0pt}{13pt}\arg(a-\gamma(t))-\arg(\gamma'(t))\right)\;,
\end{eqnarray*}
with $\rho(t):=1$ (or $\rho(t):=-1$) if $a$ is on the left (or
right) side of the oriented line $\gamma(t)+\gamma'(t)\R$. We have
a change of sign whenever $a\in\sigma_t(-\infty,t)$: in that case
$a\in \gamma(t)+\gamma'(t)\R$, and the derivative vanishes.

Differentiating under the integral (which is allowed since the
integrands are differentiable and all the integrals exist) shows
if $a\not\in\gamma(t)+\gamma'(t)\R$
\[
\frac{d}{d t}\int^\infty_t |d\arg(\gamma(s)-a)|\stackrel{(*)}{=}
-\rho(t)\frac{d}{d t} \arg(\gamma(t)-a)= -\rho(t)\frac{d}{d t}
\arg(a-\gamma(t))\;.
\]
Step (*) can be seen like this: the derivative of the integral is
negative, so that the sign of the factor in front of the
parentheses is positive (negative) if $\arg(\gamma(t)-a)$ is
decreasing (increasing) in $t$, and this is exactly the case if
$a$ is on the right (left) side of $\gamma(t)+\gamma'(t)\R$. If
$a\in\gamma(t)+\gamma'(t)\R$, the value of $\rho$ is not important
for us, and we have $\arg(a-\gamma(t))-\arg(\gamma'(t))=0$.
Therefore
\begin{eqnarray*}
v'(t)&=&\frac{\rho(t)}{2\pi}\frac d{d
t}\left(\rule{0pt}{13pt}\arg(a-\gamma(t))-\arg(\gamma'(t))\right)
-\frac{\rho(t)}{2\pi}\left(\frac{d}{d t} \arg(a-\gamma(t))\right)=\\
&&=\frac{-\rho(t)}{2\pi} \left(\frac{d}{d
t}\arg(\gamma'(t))\right)
 \qquad\mbox{and}\\
u'(t)&=&-\frac1{2\pi} \left|\frac{d}{d t}\arg(\gamma'(t))\right|
\;.
\end{eqnarray*}
Hence $u$ and $v$ are continuously differentiable with $v'(t)=\pm
u'(t)$ everywhere. Since we have $v(\infty)=1/2=u(\infty)+1/2$,
this yields $v(t_1)\le u(t_1)+1/2$. Together with
$\gamma([t_1,\infty))\subset\sigma_{t_1}(\R)$, we therefore get
$\alpha(\gamma,a)\le v(t_1) \le u(t_1)+1/2 =
\alpha(\gamma',0)+1/2$. \qed
\end{proof}

\begin{proposition}[Variation Numbers and Pullbacks]
\label{windnr_pullback} \lineclear 
Let $\gamma_0:(t_0,\infty)\to\C$ be an admissible curve
such that $\alpha((\gamma_0)',0)\le 1/2$. For every $n\in\N$
choose $a_n\in\C\sm \overline{\gamma_n(t_0,\infty)}$ and
$\gamma_{n+1}:(t_0,\infty)\to\C$ such that
$\exp(\gamma_{n+1}(t))=\gamma_n(t)-a_n$. Then
\begin{eqnarray*}
&\forall n\in\N\mbox{, }\forall t_1>t_0\mbox{, }\forall b\not\in
\overline{\gamma_n(t_1,\infty)}:\\
& \alpha(\gamma_n|_{(t_1,\infty)},b)\le 2^n\quad\mbox{and}\quad
\alpha((\gamma_n)'|_{(t_1,\infty)},0)\le 2^n-1/2\ .
\end{eqnarray*}
\end{proposition}

\begin{proof}
It follows from the calculations below that all the variation
numbers are indeed defined. The case $n=0$ follows from Lemma
\ref{BoundVar}. Let $a$, $\gamma$, $\tilde\gamma$ denote $a_n$,
$\gamma_n$, $\gamma_{n+1}=\log(\gamma_n-a)$ respectively. We
estimate
\begin{eqnarray*}
\alpha({\tilde
\gamma}',0)&=&\alpha\left(\frac{\gamma'}{\gamma-a},0\right)=
\frac1{2\pi}\int\left| d\left(\Im\log\left(\frac{\gamma'}{\gamma-a}\right)\right)\right|=\\
&=& \frac1{2\pi}\int\left|\rule{0pt}{13pt}d(\Im\log \gamma'-\Im\log(\gamma-a))\right|= \\
&=& \frac1{2\pi}\int\left|\Im\left(\frac{\gamma''}{\gamma'}
\right)-\Im\left(\frac{\gamma'}{\gamma-a}
\right)dt\right| \le \\
&\le& \frac1{2\pi}\int\left|\Im\left(\frac{\gamma''}{\gamma'}
\right)dt\right|
+\frac1{2\pi}\int\left|\Im\left(\frac{\gamma'}{\gamma-a} \right)dt\right|=\\
&=& \alpha(\gamma',0)+\alpha(\gamma,a)\stackrel{\mbox{\scriptsize
Lemma }\ref{BoundVar}}{\le}
2\alpha(\gamma',0)+1/2\le\\
&\le& 2(2^n-1/2)+1/2=2^{n+1}-1/2\ .
\end{eqnarray*}
Lemma \ref{BoundVar} thus gives us
$\alpha(\gamma_n,a)\le\alpha((\gamma_n)',0)+1/2\le 2^n$ for all
$n\in\N$. \qed
\end{proof}

\hide{
Recall that, for fixed $\s\in\S_0$, by Proposition
\ref{d2dtdynray} there is a function $(K,t)\mapsto n_{\s}(K,t)$
such that for all $|\kappa|\le K$, $n\ge n_{\s}(K,t_0)$, $t\ge
t_0>t_{\s}$
\begin{equation}\label{d2dtest'}
\left|\frac{g''_{\kappa,\sigma^{n}\s}(F^{\circ n}(t))}
{g'_{\kappa,\sigma^{n}\s}(F^{\circ n}(t))}\right|<e^{-F^{\circ
n}(t)/2}\;.
\end{equation}
This function is monotonically increasing in $K$ and decreasing in
$t$.
}

\subsection{The Bound on Parameter Rays}

In this final subsection, we will complete the proof of
Proposition~\ref{parabound}. It may be helpful to outline the general line of
argument before going into details. We want to construct parameter rays
$G_{\s}\colon(t_{\s},\infty)\to\C$; we know these exist as curves for sufficiently
large potentials (Proposition~\ref{parrayend}). The danger is that there is a
$\tilde t>t_{\s}$ such that as $t\searrow\tilde t$, $G_{\s}(t)\to\infty$.

We will of course use our estimates on dynamic rays (Theorem~\ref{dynrays}). The
problem is that these estimates depend on $\kappa$ and become worse as
$\kappa\to\infty$. Rescue comes from $\infty$ in a different way: for a given
parameter $\kappa=G_{\s}(t)$ with $t>t_{\s}$, one needs to iterate the dynamic ray
$g_{\s}(t,\infty)$ only a finite number of times until the iterated image ray is
almost horizontal: if $N$ is this number of iterations, then the winding number of
$g_{\s}(t,\infty)$ will be bounded by $2^N$ (Proposition~\ref{windnr_dynrays}).
This induces a partition of the dynamical plane with horizontal uncertainty of
approximately $2^N$ (Lemma~\ref{bound_critorb}). But if now $\Re(\kappa)$ is too
large, then the imaginary bounds imply that the singular orbit must escape very
fast (Proposition~\ref{behavsing}), which means it must have large potential $t$.
For bounded potential $t$, this yields an upper bound for
$\Re(\kappa)$. 

The success of this argument depends on the fact that as $|\kappa|$ increases, the
number $N$ of necessary iterations grows extremely slowly in $|\kappa|$, much
slower than the errors arising from the growth of $|\kappa|$ itself.

\begin{proposition}[Variation Number of Dynamic Rays]
\label{windnr_dynrays} \lineclear 
Consider an arbitrary parameter $\kappa$ and an
external address $\s\in\S_0$. Let $z=g_{\s}^{\kappa}(t_0)$ be any
point on the dynamic ray of address $\s$ with potential
$t_0>t_{\s}$. 
Suppose $N\in\N$ is such that for all $t\ge F^{\circ N}(t_0)$,
\[
\left|\frac{g''_{\kappa,\sigma^{N}\s}(t)}
{g'_{\kappa,\sigma^{N}\s}(t)}\right| 
<  
e^{-t/2}
\;.
\]
Then 
\[
\alpha\left(g_{\s}^{\kappa}|_{(t_0,\infty)},z\right)\le 2^N\ .
\]
\end{proposition}

\begin{proof}
On ray tails, we have bounds on the rays (Theorem \ref{dynrays} (\ref{defined})),
their first derivatives (Proposition~\ref{ddtdynray}) and their second derivatives
(Proposition~\ref{d2dtdynray}), so ray tails (and hence entire dynamic rays) are
admissible curves. 

Specifically, the curve 
$\gamma:=g_{\sigma^{N}\s}^{\kappa}:(F^{\circ N}(t_0),\infty)\to\C$ satisfies
\begin{eqnarray*}
\alpha(\gamma',0)
&=&
\frac1{2\pi}\int_{F^{\circ n}(t_0)}^\infty
\left|\Im\left(\frac{\gamma''(t)}{\gamma'(t)}\right)\right|dt \le
\frac1{2\pi}\int_{F^{\circ n}(t_0)}^\infty
\left|\frac{g''_{\kappa,\sigma^{n}\s}(t)}
{g'_{\kappa,\sigma^{n}\s}(t)}\right|dt < 
\\
&<& 
\frac1{2\pi}\int_{F^{\circ n}(t_0)}^\infty e^{-t/2} dt\le
\pi^{-1} e^{-F^{\circ n}(t_0)/2} < 1/2 \,\,.
\end{eqnarray*}
If we define $\gamma_0:=\gamma$ and
$\gamma_{k+1}:=L_{\kappa,s_{n-k}}(\gamma_k)$, then
$$\gamma_n=g_{\s}^{\kappa}\left.\left(F^{\circ (-n)}(t)\right)\right|_{t>F^{\circ
n}(t_0)}\;,$$and applying Proposition \ref{windnr_pullback}
settles the claim. \qed
\end{proof}

\begin{lemma}[Bounding the Imaginary Parts]
\label{bound_critorb} \lineclear 
Suppose $s\in\S_0$ and suppose $\kappa$ is a parameter
such that $g_{\s}^{\kappa}(t_0)=0$ for some potential
$t_0>t_{\s}$. Let $N$ be as in Proposition~\ref{windnr_dynrays}. 
Then the singular orbit $(z_k)_{k\ge 1}=(E_{\kappa}^{\circ (k-1)}(0))_{k\ge 1}$
satisfies for all $k\ge 1$
\begin{equation}\label{Imzk_estim}
|\Im(z_k)|\le 2\pi(2^N+1+|s_k|)+|\Im\kappa|\;.
\end{equation}
Furthermore,
\[
|\Im\kappa|\le 2\pi(2^N+1+|s_1|)\;.
\]
\end{lemma}

\begin{proof}
The inverse images of the ray $g_{\s}^{\kappa}|_{(t_0,\infty)}$
provide a \emph{dynamic partition} of the plane, as opposed to the
\emph{static partition} introduced in Section \ref{sec_dynrays}.
Note that dynamic rays cannot intersect the boundary of the
dynamic partition, so that each ray $g_{\s'}^{\kappa}$ has to be
contained in one of the two components which are asymptotic to the line
$t-\kappa+2\pi i s_1'$ for large $t$. By Proposition
\ref{windnr_dynrays}, the vertical variation of any boundary component 
of the dynamic partition is bounded by $2\pi\cdot 2^N$; since the
boundaries are a vertical distance $2\pi$ apart, (\ref{Imzk_estim})
follows.

The additional inequality follows similarly: the strip of the
dynamic partition containing $0$ is asymptotic to the line
$t-\kappa+2\pi is_1$, so that the bound on the vertical variation
within the strip yields $|-\Im\kappa+2\pi s_1|\le 2\pi (2^N+1)$,
and the triangle inequality gives the desired estimate. \qed
\end{proof}

\begin{proposition}[The Behavior of the Singular Orbit]\label{behavsing}
\lineclear Let $\kappa\in I$ be a parameter such that
$\Re\kappa\ge 4$ and $|\Im\kappa|\le \Re \kappa-2$. Let
$K:=|\kappa|$ and let $(z_k)_{k\ge
1}:=(E_{\kappa}^{\circ(k-1)}(0))_{k\ge 1}$ be the singular orbit
and suppose that for every $k\ge 2$
\begin{eqnarray*}
|\Im(z_k)| &<& F^{\circ(k-1)}(\Re\kappa-2)-K\;. \\
\mbox{Then}\quad \forall k\ge 1 \quad |\Re(z_k)|&\ge& F^{\circ
(k-1)}(\Re\kappa-1)-K\;. 
\end{eqnarray*}
\end{proposition}

\begin{proof}
Let $(B_k), (C_k)$ denote the statements
\begin{eqnarray*}
|\Re(z_k)|\ge  F^{\circ (k-1)}(\Re\kappa-1)-K && (B_k)\;,\\
\Re(z_k)\ge -\Re\kappa  && (C_k)\;.
\end{eqnarray*}
The induction seeds $(B_1):\ 0\ge \Re\kappa -1-K$ and $(C_1): 0\ge
-\Re\kappa$ are trivial. The induction steps follow immediately
from \cite{rs1}, Lemmas 4.2 and 4.3, which say, translated (and
weakened) from the $\exp(z)+\kappa$ into the $\exp(z+\kappa)$
parametrization: if $\kappa$ is a parameter as specified above,
and if there is a $k\ge 2$ such that $(C_{k'})$ holds for all $1\le
k'\le k-1$, then $(B_k)$ and $(C_k)$ hold (since $\kappa$ does not
admit an attracting orbit).\qed
\end{proof}

\remark
While the complete proof of the preceding proposition needs the detailed arguments
from \cite{rs1}, the idea is simple: if the real part of $z_k$ is large and
positive, then $|z_{k+1}|=|E_\kappa(z_k)|$ is exponentially large. If the
imaginary parts of the orbit are bounded, then $|\Re(z_{k+1})|$ must be almost as
large as $|z_{k+1}|$. If the real part is negative, then $z_{k+2}$ is extremely
close to the origin, and there is an attracting orbit of period at most $k+2$.


\hide{
\begin{proposition}[Bound on Singular Value]
\label{Prop:BoundSingValue} \lineclear
Suppose that $\kappa\in \G_{\s}(t)$ with $\s\in \S_0$ and $t>t_\s$. Then
$|\kappa|<R$ for a constant $R$ depending only on $\s$ and $t$; for sufficiently
large $t$, we have $|\kappa|<2t$.
\end{proposition}
\proof
}

\proofof{Proposition \ref{parabound}}
It suffices to consider the case $|\kappa| > e$. Set $K:=|\kappa|$. By
Proposition~\ref{d2dtdynray}, there are universal constants $A,B,C\geq 1$ such
that 
\[
\left|\frac{(g_{\sigma^n(\s)}^\kappa)''(t')}{(g_{\sigma^n(\s)}^\kappa)'(t')}\right|
<e^{-t'/2}<1
\]
for all $t'\geq t$ provided $t'\ge At_{\sigma^n(\s)}^*+B\log K+C$.

There is an $N_0\in\N$ so that $F^{\circ N_0}(t)\ge At_{\sigma^{N_0}(\s)}^*+C+1$.
Let $N_1\in\N$ be minimal with $F^{\circ N_1}(1)\ge B\log K>1$. Then $N_1\ge 1$
and $F^{\circ (N_1-1)}(1) < B\log K$.

By Lemma~\ref{Lem:BoundInitialSteps}, we have for all $n\ge N_0+N_1=: N$
\[
F^{\circ n}(t)\ge At_{\sigma^n(\s)}^*+B\log K+C \;.
\]
There is a constant $c_1>0$ such that all $n\in\N$ satisfy $2^n\le c_1 F^{\circ
(n-1)}(1)$. 

By Proposition~\ref{bound_critorb}, we have the estimates
\begin{eqnarray*}
|\Im\kappa|
&\le& 
2\pi(2^N+1+|s_1|)
\le
2\pi(2^{N_0}c_1 F^{\circ(N_1-1)}(1)+1+|s_1|)
\\
&\le& 2\pi(2^{N_0}c_1B\log K+1+|s_1|)
= c_2+c_32^{N_0}\log K
\end{eqnarray*}
with constants $c_2,c_3>0$ depending only on $|s_1|$, and also
\begin{equation}
|\Im z_k|\le 2\pi(2^N+1+|s_k|)+|\Im\kappa| \;.
\label{Eq:ImagPartsOrbit_z_k}
\end{equation}
There is an $M>0$ so that if $K\ge M$, then
\begin{eqnarray}
K &\ge& 2c_2+2c_32^{N_0}\log K+2
\label{Eq:K-bound 1}
\\
K & \ge& 
2\pi t_{\s}^*+(2\pi2^{N_0}c_1B+c_32^{N_0}+1)\log K+2\pi+4+c_2
\label{Eq:K-bound 2}
\;.
\end{eqnarray}
Now suppose $K\ge M$. This implies
\[
|\Re\kappa|\ge K - |\Im\kappa|
\ge K-c_2-c_32^{N_0}\log K
\ge c_2+c_32^{N_0}\log K+2
\ge |\Im\kappa|+2
\]
and
\begin{eqnarray*}
|\Re(\kappa)|-2
&\ge&
K - |\Im(\kappa)| -2
\\
&\ge&
2\pi t_{\s}^*+2\pi2^{N_0}c_1B\log K + \log K+2\pi+ 2
\\
&\ge&
2\pi t_{\s}^*+2\pi2^{N_0}c_1 F^{\circ(N_1-1)}(1) +2\pi+ F^{-1}(2K)
\\
&\ge&
2\pi t_{\s}^*+2\pi2^{N_0}2^{N_1} +2\pi+ F^{-1}(2K).
\end{eqnarray*}
Since parameters with $\Re(\kappa)<-1$ are known to be attracting, we conclude
that $\Re(\kappa)\ge 2\pi+2$. We obtain for $k\geq 2$ (using convexity of $F$)
\begin{eqnarray*}
F^{\circ(k-1)}(\Re\kappa-2)
&\ge&
F^{\circ(k-1)}(2\pi t_\s^*+2\pi 2^N+2\pi+F^{-1}(2K))
\\
& > &
2\pi F^{\circ(k-1)}(t_\s^*)+2\pi 2^N+2\pi+2K
\\
&\ge&
2\pi t_{\sigma^{k-1}(\s)}^*+2\pi 2^N+2\pi+2K
\\
&\ge&
2\pi(|s_k|+ 2^N+1)+2K
\ge
|\Im z_k|+K
\;.
\end{eqnarray*}

By Proposition~\ref{behavsing}, it follows for all $k\geq 2$
\[
F^{\circ(k-1)}(\Re\kappa-1)-K\le\Re z_k\stackrel{\mbox{\scriptsize
(\ref{lowpotasymp}) in Thm \ref{dynrays}}}{\le}
F^{\circ(k-1)}(t)-\Re \kappa+O(e^{-F^{\circ(k-1)}(t)})\;.
\]
Comparing the growth of the left and the right hand sides as
$k\to\infty$, we conclude $\Re \kappa\le t+1$. The triangle inequality thus yields
$K\le |\Re\kappa|+|\Im\kappa|< 2|\Re\kappa|-2\le 2t$.

Every fixed choice of $t>t_\s$ yields a fixed value of $N_0$ and thus a fixed
value of $M$; and clearly we can choose $M$ so that it depends continuously on
$t$. Then
\[
K\le\max\{2t,M\}
\hide{
K\le\max\left\{2t,2\pi t_{\s}^*+(2\pi2^{N_0}c_1B+c_32^{N_0}+1)\log K+2\pi+4+2c_2
\right\}
}
\;.
\]
Note that as $t$ increases, $N_0$ and hence $M$ do not increase, while
$c_1,c_2,c_3$ are independent of $t$. Therefore, for large $t$ we have the bound
$K\le 2t$.
\qed

\end{document}